\documentclass[pdflatex,sn-mathphys-num]{sn-jnl}
\usepackage{alphabeta}
\usepackage{graphics}
\usepackage{graphicx} 
\usepackage{rotating}
\usepackage{dcolumn}
\usepackage{hyperref}
\usepackage{longtable}
\usepackage{multirow}
\usepackage{amsmath}
\usepackage{amssymb}
\usepackage{xspace}
\usepackage{textcase}
\usepackage{float}
\usepackage{amsthm}

\usepackage{standalone}
\usepackage{tikz}
\usepackage{tikzscale}
\usepackage{pgfplots}
\pgfplotsset{compat=1.18}
\usetikzlibrary{trees}
\usetikzlibrary{arrows.meta}
\usetikzlibrary{positioning}
\usetikzlibrary{arrows, math, angles, quotes, shapes.geometric, calc}
\usepackage{commath}
\usepackage{xcolor}
\usepackage{hyperref}
\usepackage{subfiles}

\usepackage{lineno}

\newcommand{\Rbb}{\mathbb{R}}
\newcommand{\bbR}{\Rbb}
\newcommand{\Cbb}{\mathbb{C}}
\newcommand{\bbC}{\Cbb}

\newcommand{\KO}{\mathcal{K}}

\newcommand{\from}{\colon}
\newcommand{\map}[3]{#1\from #2 \to #3}
\newcommand{\mapself}[2]{\map{#1}{#2}{#2}}

\newcommand{\Linf}{\Lspace{\infty}}
\newcommand{\Lspace}[1]{\mathcal{L}^#1}

\theoremstyle{definition}
\newtheorem{definition}{Definition}[section]

\theoremstyle{thmstyleone}%
\theoremstyle{thmstyletwo}%

\theoremstyle{thmstylethree}%

\begin{document}

\title[Operator-Theoretic Methods for Differential Games]{Operator-Theoretic Methods for Differential Games} 


\author*[1]{\fnm{Craig} \sur{Bakker}}\email{craig.bakker@pnnl.gov}

\author[1]{\fnm{Adam} \sur{Rupe}}\email{adam.rupe@pnnl.gov}

\author[2]{\fnm{Alexander} \sur{Von Moll}}\email{alexander.von\_moll@us.af.mil}
\author[2]{\fnm{Adam R.} \sur{Gerlach}}\email{adam.gerlach.1@afrl.af.mil}

\affil*[1]{\orgdiv{National Security Directorate}, \orgname{Pacific Northwest National Laboratory}, \orgaddress{\street{902 Battelle Boulevard}, \city{Richland}, \postcode{99352}, \state{Washington}, \country{United States of America}}}

\affil[2]{\orgdiv{\orgname{Air Force Research Laboratory}, \orgaddress{\city{Wright-Patterson Air Force Base}, \postcode{45433}, \state{Ohio}, \country{United States of America}}}}

\abstract{
Differential game theory offers an approach for modeling interactions between two or more agents that occur in continuous time. The goal of each agent is to optimize its objective cost functional. In this paper, we present two different methods, based on the Koopman Operator (KO), to solve a zero-sum differential game. The first approach uses the resolvent of the KO to calculate a continuous-time global feedback solution over the entire domain. The second approach uses a discrete-time, data-driven KO representation with control to calculate open-loop control policies one trajectory at a time. We demonstrate these methods on a turret defense game from the literature, and we find that the methods' solutions replicate the behavior of the analytical solution provided in the literature.. Following that demonstration, we highlight the relative advantages and disadvantages of each method and discuss potential future work for this line of research.}


\maketitle

\section{Introduction}

Differential game theory offers an approach for modeling interactions between two or more agents that occur in continuous time; this can be used, for example, to develop optimal and robust control inputs in adversarial scenarios. In differential games, the agents make decisions continuously and simultaneously unlike one-shot games such as Rock-Paper Scissors, wherein players make one decision \textit{a priori}, or turn-based games such as Chess. In a differential game, each player has a particular objective cost functional that may depend on the state of the system, its control signals, and the control signals of other agents. The goal of each agent is to optimize (i.e., maximize or minimize, depending on the scenario) its objective cost functional. The Nash equilibrium concept states that no single player has an incentive to deviate from its equilibrium strategy when all others play their equilibrium strategies. The value of the objective cost functional when all players employ their equilibrium strategies is called the Value of the game.

In his seminal work on differential games, Rufus Isaacs developed the theory and applied it to scenarios common in warfare such as pursuit-evasion and attrition-based conflicts~\cite{isaacs1965differential}. Much of his work was devoted to a particular class of game, namely, two-player zero-sum differential games. For this class of problem, the Nash equilibrium specializes to the saddle-point equilibrium. Under a saddle-point equilibrium, a player is guaranteed to perform as good or better than the Value of the game \textit{regardless of their opponent's actions}. In other words, an agent who plays according to its equilibrium strategy is \textit{robust} with respect to its adversary's strategies. This robustness is what makes differential game solutions (to two-player, zero-sum scenarios) so desirable.

Although differential game theory has been applied extensively in the field of economics~\cite{dockner2000differential}, most of the applications within engineering have been to pursuit-evasion type scenarios. The Lady in the Lake scenario, for example, involves an agent who begins inside a circle, has single integrator kinematics, and seeks to reach the circle with maximum angular separation from its opponent who is constrained to move along the circle~\cite{basar1982chapter,gardner1975mathematical}. The latter agent seeks to minimize the angular separation at the time of exit of the former. Another famous example is the Homicidal Chauffeur scenario wherein a faster, turn-constrained agent seeks to intercept a slower, agile agent; the agents seek to minimize or maximize, respectively, the time to intercept~\cite{merz1971homicidal}. More recently, differential game theory has been applied to engagements between stationary, turn-constrained \textit{turrets} and mobile attackers~\cite{akilan2017zero-sum}. Much of the existing literature on differential games has focused on games with relatively simple dynamic models and few parameters in order to make the analysis tractable. This is in stark contrast with recent learning-based approaches for obtaining strategies in multi-agent adversarial settings which have been used to obtain strategies for much more complex scenarios~\cite{zhang2021multi-agent}.

When it comes to solving differential games, there are two significant branches in the literature: traditional methods that are capable of addressing nonlinear systems and methods for linear systems. The former methods are powerful in the sense that their solutions yield the equilibrium control action for each player and the resulting equilibrium cost/reward value of the game. However, these methods suffer from the curse of dimensionality and do not scale well for systems of even modest fidelity~\cite{bernhard2014pursuit-evasion}.
For example, the Homicidal Chauffeur Differential Game involves an agents with first- and second-order kinematics, and it took 17 years to obtain a full solution~\cite{merz1971homicidal}. The main difficulty in solving differential games with nonlinear dynamics lies in properly addressing the various singular surfaces that arise~\cite{basar1982chapter, isaacs1969differential}. Systems with more complex dynamics typically have more singularities. At present, there appears to be no general way to automate the process of solving a differential game with singularities. Furthermore, unless the solution admits a closed-form, state-feedback solution (which becomes less likely as model fidelity increases) implementing the control in real time becomes essentially intractable.  

On the other hand, if one can obtain a linear representation of the system dynamics, then more standard linear game-theoretic tools are available (e.g., LQ games). Those linear game-theoretic solutions, moreover, may be in a state-feedback closed form~\cite{basar1982chapter6}, which makes real-time implementation feasible.  

If the linear dynamics are discrete-time rather than continuous-time, then it is sometimes also possible to solve for the game's equilibrium by transforming the game into a Mixed Complementarity Problem (MCP)~\cite{quigley2020multiagent}. If each player's optimization is convex (assuming other players' variables are fixed), the overall solution can be defined by combining the optimality conditions of each player's optimization together. The PATH algorithm~\cite{dirkse1995path} can then solve the resulting MCP directly for the problem's Nash equilibrium~\cite{ruiz2014tutorial}. This approach can be extended to different game theoretic situations via the use of Mathematical Programs with Equilibrium Constraints (MPECs) and Equilibrium Problems with Equilibrium Constraints (EPECs), though MPECs and EPECs require different solution methods~\cite{ruiz2014tutorial}.

Linear methods can easily incorporate more states/dimensions, but the assumption of linearity severely restricts the set of problems to which they can be applied (see, e.g.,~\cite{garcia2020defense}). However, the theory of the Koopman operator (KO)~\cite{koopman_hamiltonian_1931} provides a framework for reasoning about finite dimensional nonlinear systems in a linear, but infinite-dimensional, setting. The KO has gained much attention in the dynamics and controls communities as it provides a theoretical means for bridging the gap between nonlinear systems and linear solution approaches~\cite{budisic12jsr, mauroy_koopman_2020}. 

The linearity of the KO, and its control representation~\cite{kaiser_data-driven_2021}, make it potentially useful for games with embedded dynamics (e.g., differential games). To the best of our knowledge, this has only been done in a few cases. Krolicki et al. formulate robust control as a min-max problem with a continuous-time KO representation, and they solve it using policy iteration~\cite{krolicki2022koopman}. Zhao and Zhu use the KO in the context of a Stackelberg game\footnote{A sequential game form where a ``leader" acts first and ``followers" then act in response to the leader.}, but they do not actually use the KO to solve the Stackelberg game directly~\cite{zhao2023stackelberg}. Instead, they transform the Stackelberg game into a constrained optimization in its original nonlinear formulation and then apply the KO to the resulting dynamics; given that the optimality conditions that the authors use to transform the game into a constrained optimization are not sufficient in a general nonlinear context, the authors' approach raises some questions about its generalizability. By contrast, in Oster et al.~\cite{oster2023multi}, the KO linearity enables the authors to use strong duality to convert a tri-level optimization into a bi-level optimization that can then be solved using existing techniques from the multi-level optimization literature. The authors discuss the use of the Karush-Kuhn-Tucker (KKT), or optimality, conditions instead of duality but do not actually present results using the KKT conditions. Unlike Krolicki et al., Oster et al. actually consider (more than two) separate decision-making agents in their formulation.

\subsection{Contributions}

Motivated by the limitations inherent to current methods for solving differential games, we present here approaches for automating the process of solving a wide variety of games by using the KO.  Specifically, we develop two different methods, based on the KO, to solve a zero-sum differential game and demonstrate these methods on an example from the literature. 

The first approach uses the resolvent of the KO to calculate a continuous-time global feedback solution over the entire domain. To do this, we use a fixed, finite-dimensional set of basis functions to represent the Koopman generator and our feedback control policies. We calculate the action of the KO, via the resolvent, on an objective function with terminal and running costs and use sequential optimizations to calculate the game's saddle point. 

The second approach uses a discrete-time, data-driven KO representation with control that allows us to represent a) the objective as a quadratic function of the states and control inputs and b) the dynamics as linear in the lifted states (i.e., observable functions) and control inputs. We then produce an MCP from the optimality conditions for the game's players and solve the resulting MCP to calculate open-loop control policies one trajectory at a time. These two approaches represent fundamentally different ways of using the KO to solve a differential game, but the combination of them in a single paper allows us to compare and contrast them. They could be addressed individually in separate papers, but considering them together here highlights their complementary strengths and weaknesses.

Following the presentation of these methods, we demonstrate them on a zero-sum differential game drawn from the literature and assess the quality of these solutions against a ground-truth equilibrium policy.  As we show, our approaches yield the equilibrium solutions for the \textit{approximated} dynamics (i.e., using a finite approximation of the Koopman operator or generator) and thus stand in the gap between purely game-theoretic methods and purely data-driven methods.  We conclude the paper with a discussion of method computational costs, method performance, and next steps for this line of research.

The remainder of the paper is organized as follows.
Section~\ref{sec:background} provides background on the KO and MCPs. Section~\ref{sec:methodology} presents the details of the two KO-based differential game solution methods. Section~\ref{sec:turret-defense-problem} introduces the particular differential game, examined previously in the literature, with which we demonstrate our KO-based methods and describes how those methods are applied. Section~\ref{sec:results} summarizes the previously existing solution and presents the results produced by the two solution methods. Section~\ref{sec:discussion} compares the two numerical approaches, discusses some implementation details, and outlines future work.  Section~\ref{sec:conclusions} concludes the paper.

\section{Background}
\label{sec:background}

\subsection{The Koopman Operator}

\subsubsection{Definition and Properties}
\label{sec:KO defn}

Given a dynamical system, its corresponding KO is a linear operator on the space of functions that map the system states to a scalar value. This operator governs how these functions evolve by the underlying system. Formally, for discrete-time dynamics $\mapself{F}{\bbR^n}$, the KO is defined as 
\begin{equation}
    \KO_F g = g \circ F, \quad \forall g \in \Linf,~\map{g}{\bbR^n}{\bbR}, \label{eq:koopman}
\end{equation}
where $\mapself{\KO_F}{\Linf}$ is the KO corresponding to dynamics $F$ and $a \circ b$ denotes the composition of functions $a$ and $b$ such that $\del{a \circ b}\del{x} = a\del{b\del{x}}$~\cite{lasota_chaos_1994}. This operator possesses two important properties that are leveraged in this work~\cite{lasota_chaos_1994}:

\begin{enumerate}
    \item $\KO$ is a linear operator: $\KO(\alpha g + \beta h) = \alpha \KO g + \beta \KO h, \quad \forall g,h \in \Linf,~\alpha,\beta\in\bbR$
    \item $\KO$ is a \emph{contracting} operator: $\norm{\KO g}_\infty \leq \norm{g}_\infty, \quad \forall g \in \Linf$
\end{enumerate}

For continuous-time dynamics defined implicitly by the ordinary differential equation $\dot x = f(x)$ with initial condition $x(0) = x_0$, one can consider the set of possible discrete-time dynamics parameterized by the time step $\cbr{F_t}_{t\geq 0}$ where
\begin{equation}
    F_t\del{x_0} = x\del{t} = x_0 + \int_0^t f(x(\tau)) \dif{\tau} \label{eq:flowmap}
\end{equation}
is the flow function mapping $x_0$ forward in time by a duration $t$. Thus, for each fixed $t \geq 0$, one can uniquely define the KO $\mapself{\KO_t}{\Linf}$ corresponding to the discrete-time dynamical system $F_t \in \cbr{F_t}_{t\geq 0}$. This leads to the family, or \emph{semigroup}, of KOs $\cbr{\KO_t}_{t\geq 0}$.

In general, the flow map $F_t \left(x\right)$ is not known in closed form for the dynamics $\dot x = f(x)$. In these cases, it is useful to consider the so-called \emph{infinitesimal Koopman Operator}, also referred to as the \emph{Koopman Generator} (KG). The KG, $L$, is defined such that 
\begin{equation}
    \od{}{t} \sbr{\KO_t g} = L \sbr{\KO_t g}, \quad \KO_0 g = g
\end{equation}
where 
\begin{equation}
    L = \lim_{t\rightarrow 0} \frac{\KO_t - I}{t} = f \cdot \nabla_x. \label{eq:koop_gen}
\end{equation}

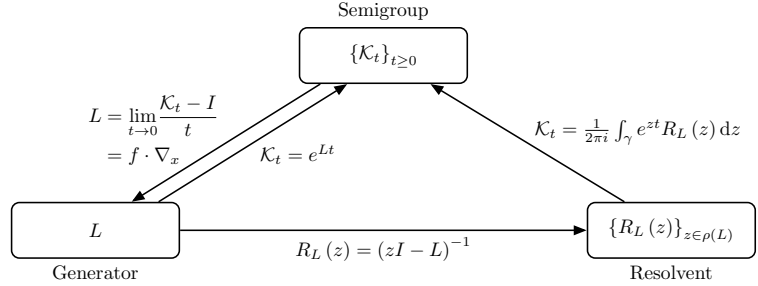
\begin{figure}[htp]
\centering
\begin{tikzpicture}

\def\yshift{0.2}
\def\nodespace{3cm}
\def\generator{L}
\def\semigroup{\mathcal{K}}

\tikzstyle{block} = [rectangle, rounded corners, thick, minimum width=3cm, minimum height=1cm, text centered, draw=black]
\tikzstyle{arrow} = [thick, ->, >={Stealth[scale=1.0, inset=0pt]}]
\tikzset{double arrows/.style args={#1-#2}{%
           insert path={(#1) edge[bend right=7,draw=none] coordinate[at start](#1-b) coordinate[at end](#2-b) (#2)
                             edge[bend left=7,draw=none] coordinate[at start](#1-t) coordinate[at end](#2-t) (#2)
 }}}

\node[label=above:Semigroup] (semigroup) [block] {$\cbr{\semigroup_t}_{t\geq 0}$};
\node[label=below:Generator] (generator) [block] [below left=\nodespace of semigroup] {$\generator$};
\node[label=below:Resolvent] (resolvent) [block] [below right=\nodespace of semigroup] {$\cbr{R_\generator\del{z}}_{z\in\rho\del{\generator}}$};

\draw[double arrows={semigroup-generator}, arrow] (generator-t) -- node[anchor=west, shift={(.0,-\yshift)}]{$\semigroup_t = e^{\generator t}$} (semigroup-t);
\draw[arrow] (semigroup-b) -- node[anchor=east,shift={(-.0,\yshift)}]{$\begin{aligned}
\generator &= \displaystyle{\lim_{t\to 0}} \frac{\semigroup_t  - I}{t}\\ &= f\cdot \nabla_x\end{aligned}$} (generator-b);

\draw [arrow] (generator) -- node[anchor=north] {$R_\generator\del{z} = \del{zI - \generator}^{-1}$}(resolvent);
\draw [arrow] (resolvent) -- node[anchor=west, shift={(0,\yshift)}] {$\semigroup_t = \frac{1}{2\pi i} \int_\gamma e^{zt} R_\generator\del{z}\dif z$} (semigroup);

\end{tikzpicture}
\caption{Relationships between the Semigroup, Generator, and Resolvent (Adapted from~\cite{engel_short_2006})}
\label{fig:operator_relations}
\end{figure}

One can recover the semigroup of KOs $\cbr{\KO_t}_{t\geq 0}$ from the KG by leveraging the linearity of the KO either via the operator exponential~\cite{engel_short_2006}
\begin{equation}
    \KO_t = e^{Lt}
\end{equation}
\noindent which can be thought of as an infinite dimensional analog to the matrix exponential, or via the contour integral
\begin{equation}
    \KO_t = \frac{1}{2\pi i} \int_\gamma e^{zt} R_L\del{z}\dif{z}\label{eq:contoureq}
\end{equation}
\noindent where 
\begin{equation}
    R_L\del{z} = \del{zI - L}^{-1}, \quad z \in \rho\del{L}
\end{equation}
\noindent is the \emph{Resolvent Operator} (RO), $\rho\del{L} = \bbC \backslash\sigma\del{L}$ is the resolvent set, $\sigma\del{L}$ is the spectrum of the KG, and $\gamma$ is a contour that encircles the spectrum of $L$ counterclockwise. The spectrum of $L$ is where invertibility of $zI  -L$ fails. This approach can be thought of as an infinite-dimensional analog to taking the Laplace and inverse Laplace transforms~\cite{horning_family_2024}. The relationships between the KO semigroup, KG, and RO are illustrated in Figure~\ref{fig:operator_relations}.

Although seemingly more complex, the RO-based approach offers several numerical advantages over the exponential operator~\cite{horning_family_2024,engel_one-parameter_2000}. However, special care must be taken when applying Eq.~\ref{eq:contoureq} due to the fact that $L$ is an unbounded operator. Fortunately, the contraction property of the KO allows us to directly leverage high-order quadrature schemes with explicit error bounds for solving equations of this form~\cite{horning_family_2024}. These schemes are leveraged in this work and take the following form:
\begin{equation}
    \KO_t g \approx (2\delta - L)^m \sbr{\frac{h}{2\pi}\sum_{k=-N}^N\frac{e^{\del{\delta + ihk}t}}{\del{\delta - ihk}^m} R_L\del{\delta + ihk}}g
    \label{eq:resolv approx}
\end{equation}
where $\delta,~m,~h,~N$ are quadrature parameters relating to the location of the contour, the regularity of $g$, the discretization of the contour, and the truncation of the contour~\cite{horning_family_2024}. For a more thorough treatment of the KO, KG, RO, and their relationship, the interested reader is referred to~\cite{lasota_chaos_1994, horning_family_2024, engel_one-parameter_2000}.

\subsubsection{Finite-Dimensional Koopman Approximations for Control}

In principle, $\mathcal{K}_t$ is an infinite-dimensional object. If we assign an infinite-dimensional basis to the function space, the KO can be represented as an infinite-dimensional matrix. However, the KO may have finite-dimensional invariant subspaces of that function space: an invariant subspace $S$ satisfies
\begin{align}
    \mathcal{K}_t g \left(x\right) \in S, \quad \forall \ g \in S
\end{align}

If $S$ has a finite-dimensional basis, the projection of the KO onto $S$ is a finite-dimensional matrix (in terms of the basis for $S$). In practice, numerical methods are typically used to find or use a finite-dimensional subspace that is approximately invariant~\cite{bakker2020koopman}. However, these finite approximations only possess point spectra while the original operator can have a more complex discrete and continuous spectrum~\cite{ colbrook_rigorous_2024,kato_perturbation_1995}. Furthermore, these finite approximations can struggle with spectral pollution introduced by discretization, though the RO-based ResDMD approach helps alleviate this issue~\cite{colbrook_rigorous_2024}. 

The KO is analytically useful on its own -- e.g., for characterizing global stability~\cite{mezic15cp}, quantifying system observability~\cite{surana2016linear}, and informing model reduction~\cite{budisic12jsr}. However, the KO can also be used to represent controlled system: $\dot{x} = f \left(x,\upsilon\right)$, where $\upsilon \in \mathbb{R}^m$ is the vector of control inputs. The baseline finite Koopman approximation
\begin{gather}
    \Psi \left(x_{t+1}\right) = K \Psi \left(x_t\right), \quad \Psi \in S
\end{gather}
\noindent can be extended to include control in several different ways. The fully linear form
\begin{align}
    \Psi \left(x_{t+1}\right) &= K \Psi \left(x_t\right) + K_\upsilon \upsilon_t, \quad \Psi \in S
    \label{eq:linear KO control}
\end{align}
\noindent is perhaps the most common control representation used because it keeps the dynamics linear in both the lifted states and the control inputs, but it does have other limitations~\cite{bakker2019koopman}. If there is a control objective that is linear or quadratic in the control inputs and lifted states, and if any constraints are linear in the lifted states and control inputs, then the resulting optimal control problem is a linear or quadratic program; this is likely to be much easier to solve, even if the lifting process increases the dimensionality of the problem, than the original nonlinear optimal control problem. This kind of formulation can also be modified to include domain knowledge, thereby improving the quality of the approximation and ensuring that the Koopman dynamics reproduce key features of the original dynamics~\cite{king2021solving}. Once this formulation is in place, it can be used for optimal control~\cite{king2021solving,korda2018linear}, for quantifying the resilience of the dynamical system in question~\cite{ramachandran2022data,sinha2022data}, or for system stabilization~\cite{huang2018feedback}.

This finite approximation is just that, though -- an approximation. In particular, this formulation assumes that it is possible, at least approximately, to separate out the impact of the control inputs on the lifted state from the lifted state itself. However, the KG is actually a function of the control inputs:
\begin{align}
    L &= \sum_i f_i \left(x,u\right) \frac{\partial}{\partial x_i}
\end{align}

\noindent and thus so, too, is the KO. This recognition has clear relevance to the use of the resolvent for KO calculations. In particular, it suggests that a resolvent-based control approach that accounts for the dependence of the KG on the control inputs might be able to produce better quality solutions than an approach that relies on Eq.~\ref{eq:linear KO control}, as there is one less approximating assumption being used. However, to the best of our knowledge, no one has used the resolvent directly for optimal control. Some recent work has focused on resolvent computations and analysis (e.g.,~\cite{herrmann2021data,susuki2024koopman}), but the use of the KO for optimal control has focused more on the use of the form in Eq.~\ref{eq:linear KO control} and variants thereof.

\subsection{Mixed Complementarity Problems}

MCPs are a way of formulating game theory and optimization problems. They are most commonly used in energy market modelling~\cite{huppmann2014market}, but they have also been used to model multi-agent interactions related to food security~\cite{bakker2018shocks} and counter-smuggling~\cite{bakker2020multi}, for example. They are a computationally efficient way of tackling large-scale game theory problems -- handling hundreds of thousands of variables -- where the agents' decision variables are continuous. In particular, if the agent optimizations are linear or quadratic programs, the resulting MCP is linear, which makes it quicker and easier to solve. As such, MCPs are often deliberately formulated in this fashion. Other convex, differentiable agent optimizations can be used, though.

Optimization problems are typically solved using a numerical optimizer, and there are good gradient-based methods (e.g., IPOPT~\cite{wachter2006implementation}) for solving constrained nonlinear optimization problems. These methods can be used for nonconvex problems, but convex problems like linear or quadratic programs are particularly easy to solve, as their global optima can be found using gradient-based optimizers; in non-degenerate cases, the global optimum is unique.

We can also solve optimization problems via their KKT conditions. For a general, continuous nonlinear optimization
\begin{align}
    \min_x & f\left(x\right) \\
    g\left(x\right) & \leq 0 \\
    h \left(x\right) & = 0
\end{align}
these optimality conditions are
\begin{align}
    \frac{\partial f}{\partial x} + \mu^T \frac{\partial g}{\partial x} + \lambda^T \frac{\partial h}{\partial x} &= 0 \\
    0 \leq - g\left(x\right) \perp \mu & \geq 0 \\
    h \left(x\right) &= 0
\end{align}
where $0 \leq y \perp z \geq 0$ indicates that $y\geq 0$, $z \geq 0$ and $yz = 0$. If this problem is convex, then the KKT conditions are both necessary and sufficient to identify the solution to the original optimization problem; if the problem is not convex, then the conditions may be necessary but not sufficient, i.e., there may be multiple points that satisfy the KKT conditions. This is a \emph{mixed} complementarity problem because it is a mixture of simple equality constraints and complementarity constraints. Such problems can be solved using the PATH algorithm~\cite{dirkse1995path}, but it is generally faster to solve single-agent, single-objective optimizations directly using a numerical optimizer rather than solving the MCP.  MCPs are mostly used when the problem is game theoretic -- i.e., there are multiple simultaneous, interconnected optimizations taking place.  In such cases, it is often not possible to solve for equilibria using a standard numerical optimizer.

\section{Methodology}
\label{sec:methodology}

Let us assume that we have a two-player min-max differential game:
\begin{align}
    \min_v \max_u & \ g\left(x_T\right) + \int \limits_0^T h \left(x_t\right) dt \\
    \dot{x} &= f \left(x_t,u_t,v_t\right) \\
    0 & \geq c\left(x_t,u_t,v_t\right)
\end{align}
 where $x\in X$ is the state, $u \in U$ and $v \in V$ are the control policies for the two players, and the game has a fixed time horizon $T$; let us also define the total cost functional as $J$:
\begin{align}
    J \left(x_0,T\right) &= g \left(F_T\left(x_0\right)\right) + \int \limits_0^T h \left(F_t\left(x_0\right)\right) dt \\
    J \left(T\right) &= \int_X \left[g \left(F_T\left(x\right)\right) + \int \limits_0^T h \left(F_t\left(x\right)\right) dt \right] dx
\end{align}

\noindent where $F_t$ is defined according to Eq.~\ref{eq:flowmap}. Our goal is to solve for the equilibrium control policies $u^*$ and $v^*$ that solve this game.

\subsection{Discrete-Time, Complementarity-Based Approach}
\label{sec:edmd_methodology}

Let us assume that we have a set of basis functions $\Psi$, where $\Psi = [\psi_1 \; \psi_2 \; \cdots \; \psi_N]^\top$; we call $\Psi(x) = [\psi_1(x) \; \psi_2(x) \; \cdots \; \psi_N(x)]^\top$ the \emph{lifted state} vector. We can then represent the dynamics of the game using a discrete-time Extended Dynamic Mode Decomposition with control (EDMDc) model with a lifted state vector that has an additive dependence on the unlifted controls:
\begin{align}
    \Psi \left(x_{t+1}\right) &= K \Psi \left(x_t\right) + K_u u_t + K_v v_t
\end{align}

EDMDc is the best-fit set of $K$, $K_u$, and $K_v$ matrices that evolves the lifted state vectors, where $K$ is the part of the KO approximation associated with the effect of the lifted states on their own evolution and $K_u$ and $K_v$ are the parts of the KO approximation associated with the effects of the players' control inputs on the lifted states' evolution (see Eq.~\ref{eq:linear KO control}). That best-fit solution is found through least-squares regression -- see~\ref{app:edmd} for more details. If the components of the cost function only depend on the states, as in the formulation above, then we can also define a quadratic approximation to those components in terms of the lifted states:
\begin{align}
    g \left(x_T\right) & \approx \zeta_T^T Q_g \zeta_T \\
    \int \limits_t^{t + \Delta t} h \left(F_\tau\left(x_0\right)\right) d\tau & \approx \Delta t \zeta_t^T Q_h \zeta_t \\
\end{align}
\noindent where
\begin{align}
    \zeta_t &= \left\{ \begin{array}{c} \Psi \left(x_t\right) \\ 1 \end{array} \right\} \\
    Q_g &= \left[ \begin{array}{cc} Q_{g,\psi,\psi} & Q_{g,\psi,1} \\ Q_{g,1,\psi} & Q_{g,1,1} \end{array} \right] \\
    Q_h &= \left[ \begin{array}{cc} Q_{h,\psi,\psi} & Q_{h,\psi,1} \\ Q_{h,1,\psi} & Q_{h,1,1} \end{array} \right]
\end{align}

The running and terminal cost approximations include a constant term in $\zeta$, and this allows us to create a fully generalized second-order objective that includes linear terms as well as quadratic ones; in principle, we could restrict the cost approximations to be linear if so desired. Note that, for example, $Q_{g,1,1}$ is a scalar, whereas $Q_{g,\psi,1}$ is a column vector. The differential game is therefore approximated as
\begin{align}
    \max_{u_t} \min_{v_t} & \ \zeta_T^T Q_g \zeta_T + \Delta t \sum_{t=1}^T \zeta_t^T Q_h \zeta_t \\
    \Psi_{t+1} &= K \Psi_t + K_u u_t + K_v v_t  &t = 1,\ldots, T-1 \\
    0 & \geq c \left(\Psi_t,u_t,v_t\right)
\end{align}
where $\Psi_t = \Psi \left(x_t\right)$; we assume $\Psi$ is injective. The optimality conditions are
\begin{align}
    0 =& -\Delta t \left[\left(Q_{h,\psi,\psi} + Q^T_{h,\psi,\psi}\right) \Psi_t + \left(Q_{h,\psi,1} + Q^T_{h,1,\psi}\right)\right]  \nonumber \\
    &+ \mu_t^T\frac{\partial c \left(\Psi_t,u_t,v_t\right)}{\partial \Psi} - \lambda_{u,t-1} + \lambda_{u,t}^T K  & t = 2,\ldots,T-1 \\
    0 =& -\left[\left(Q_{g,\psi,\psi} + Q^T_{g,\psi,\psi}\right) \Psi_t + \left(Q_{g,\psi,1} + Q^T_{g,1,\psi}\right)\right]  \nonumber \\
    &+ \mu_T^T \frac{\partial c \left(\Psi_T,u_T,v_T\right)}{\partial \Psi} - \lambda_{u,T-1} & t = T \\
    0 =& \Delta t \left[\left(Q_{h,\psi,\psi} + Q^T_{h,\psi,\psi}\right) \Psi_t + \left(Q_{h,\psi,1} + Q^T_{h,1,\psi}\right)\right]  \nonumber \\
    &+ \mu_t^T\frac{\partial c \left(\Psi_t,u_t,v_t\right)}{\partial \Psi} - \lambda_{v,t-1} + \lambda_{v,t}^T K  & t = 2,\ldots,T-1 \\
    0 =& \left[\left(Q_{g,\psi,\psi} + Q^T_{g,\psi,\psi}\right) \Psi_t + \left(Q_{g,\psi,1} + Q^T_{g,1,\psi}\right)\right]  \nonumber \\
    &+ \mu_T^T \frac{\partial c \left(\Psi_T,u_T,v_T\right)}{\partial \Psi} - \lambda_{v,T-1} & t = T \\
    0 \leq& -c \left(\Psi_t,u_t,v_t\right) \perp \mu_t \geq 0 & t = 2,\ldots,T &  \\
    0 =& \lambda_{u,t}^T K_u  + \mu_t^T\frac{\partial c \left(\Psi_t,u_t,v_t\right)}{\partial u} & t = 1,\ldots,T-1 \\
    0 =& \lambda_{v,t}^T K_v  + \mu_t^T\frac{\partial c \left(\Psi_t,u_t,v_t\right)}{\partial v} & t = 1,\ldots,T-1 \\
    0 =& -\Psi_{t+1} + K \Psi_t + K_u u_t + K_v v_t & t = 1,\ldots,T-1
\end{align}

Note that $\lambda_u$ is the dual variable associated with the dynamics from the $u$ player's optimization and $\lambda_v$ is the dual variable associated with the $v$ player's optimization. In general, we would expect $\lambda_u = - \lambda_v$, but there may be exceptions to this. These optimality conditions define an MCP. Ideally, we would use observables such that the constraints $c \left(\Psi_t,u_t,v_t\right)$ would be linear in $\Psi_t$, $u_t$, and $v_t$, in which case, the MCP would be linear. We can then solve this MCP using, e.g., the PATH algorithm~\cite{dirkse1995path}.

\subsection{Continuous-Time, Resolvent-Based Approach}
\label{sec:resolve_methodology}

The baseline resolvent calculations described in Section~\ref{sec:KO defn} can be used to propagate a function forward in time. This suffices for the calculations associated with the terminal cost $g \left(x_T\right)$ of our differential game. However, integrating the running cost $h\left(x_t\right)$ requires some additional steps:
\begin{align}
    g \left(F_T\left(x_0\right)\right) = \left[\KO_T g\right]\left(x_0\right) &= \frac{1}{2 \pi i} \int_\gamma e^{zT} \left(Iz - L\right)^{-1} g \left(x_0\right) dz \\
    h \left(F_t\left(x_0\right)\right) = \left[\KO_t h\right]\left(x_0\right) &=  \frac{1}{2 \pi i} \int_\gamma e^{zt} \left(Iz - L\right)^{-1} h \left(x_0\right) dz \\
    \int \limits_0^T h \left(x\left(t\right)\right) dt &= \int \limits_0^T \frac{1}{2 \pi i} \int_\gamma e^{zt} \left(Iz - L\right)^{-1} h \left(x_0\right) dz dt \nonumber \\
    &= \frac{1}{2 \pi i} \int_\gamma e^{zT} \left(Iz - L\right)^{-1}\frac{1}{z}\left[ 1 - e^{-zT}\right] h \left(x_0\right) dz \\
    J \left(x_0,T\right) & = \frac{1}{2 \pi i} \int_\gamma e^{zT} \left(Iz - L\right)^{-1} \left( g \left(x_0\right) + \frac{1}{z}\left[ 1 - e^{-zT}\right] h \left(x_0\right) \right) dz \\
    J\left(T\right) &= \int_X \left[\frac{1}{2 \pi i} \int_\gamma e^{zT} \left(Iz - L\right)^{-1} \left( g \left(x\right) + \frac{1}{z}\left[ 1 - e^{-zT}\right] h \left(x\right) \right) dz \right] dx
\end{align}

To make the necessary calculations computationally feasible, we make several approximations. Firstly, we assume that the control policies are feedback policies: $u = u \left(x\right)$ and $v = v \left(x\right)$. Secondly, we assume a fixed set of basis functions with which to approximate the dynamics, costs, and control policies (e.g., $u\left(x\right) = \sum_k u_k \phi_k \left(x\right)$, where $\phi_i\left(x\right)$ are the basis functions). This set of basis functions need not be the same basis as used for the complementarity-based approach. The combination of these two approximations transforms the task of finding (infinite-dimensional) control policies into the task of finding the equilibrium coefficients of the control policies (e.g., $u^*\left(x\right) = \sum_k u^*_k \phi_k \left(x\right)$). It also makes the (finite approximation of the) generator $L$ a function of the control coefficients $u_k$ and $v_k$. Thirdly, we approximate the resolvent calculation using the numerical calculation scheme described in Eq.~\ref{eq:resolv approx} plus a finite sampling approach to approximate the integral over the domain $X$:
\begin{align}
    J \left(x_0, T\right) & \approx \sum_j \beta_j \left(u_k,v_k\right) \left(I z_j - L\left(u_k,v_k\right)\right)^{-1} \left( g \left(x_0\right) + \frac{1}{z_j}\left[ 1 - e^{-z_jT}\right] h \left(x_0\right) \right) \\
    J \left(T\right) & \approx \sum_{i,j} \beta_j \left(u_k,v_k\right) \left(I z_j - L\left(u_k,v_k\right)\right)^{-1} \left( g \left(x_i\right) + \frac{1}{z_j}\left[ 1 - e^{-z_jT}\right] h \left(x_i\right) \right)
\end{align}

\noindent where $\beta_j \left(u_k,v_k\right)$ are the effective weights associated with the integration sampling scheme; these weights are a function of the control inputs by virtue of the fact that they incorporate the $\left(2\delta - L \left(u_k,v_k\right)\right)^m$ term from Eq.~\ref{eq:resolv approx}.  This calculation allows us to calculate our objective, with the running cost, using the existing algorithms for resolvent-based calculations~\cite{horning_family_2024}. Our differential game definition then becomes
\begin{align}
    \min_{v_k} \max_{u_k} &\sum_{i,j} \beta_j \left(u_k,v_k\right) \left(I z_j - L \left(u_k,v_k\right)\right)^{-1} \left( g \left(x_i\right) + \frac{1}{z_j}\left[ 1 - e^{-z_jT}\right] h \left(x_i\right) \right) \\
    0 &\geq \tilde{c} \left(x_i,u_k,v_k\right)
\end{align}
\noindent where $\tilde{c} \left(x_i,u_k,v_k\right)$ is the transformation of $c \left(x,u,v\right)$ into a sampled, finite-basis form. For example, if we have a constraint $u - 1 \leq 0$, the transformed version of this would be
\begin{align}
    0 &\geq \sum_k u_k \phi_k \left(x_i\right) - 1, \quad  \forall \ i
\end{align}

To solve this (finite-dimensional) min-max game, we alternately solve each player's optimization problem with a constrained nonlinear optimization solver, the other player's control policy being held fixed, using warm starts from each previous iteration. To prevent cycling, which can happen when using this approach~\cite{wibisono2022alternating}, we limit the number of iterations before switching to the other player.

\section{The Turret Defense Problem}
\label{sec:turret-defense-problem}

\subsection{Problem Description}

\begin{figure}[htpb]
    \centering
    \includegraphics[width=0.5\textwidth]{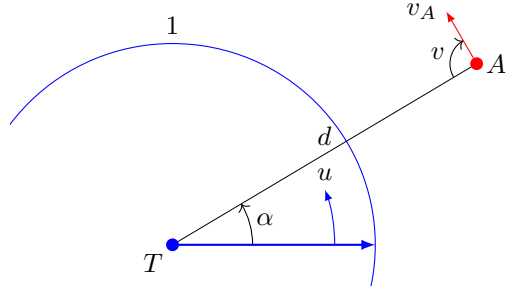}
    \caption{Turret Defense Problem Schematic}
    \label{fig:turret-defense-schematic}
\end{figure}

We begin with a modified version of the turret defense scenario investigated by Akilan et al.~\cite{akilan2017zero-sum} as presented in Von Moll et al.~\cite{vonmoll2024constrained} and shown in Figure~\ref{fig:turret-defense-schematic}:
\begin{align}
    \dot{d} &= - v_A \cos v \\
    \dot{\alpha} &= \frac{v_A}{d} \sin v - u
\end{align}
\noindent where the agent moves with constant speed, $v_A$, and controls its heading, $v$, and the turret controls its turn rate $u$. We also introduce a change of variables to avoid having a state variable in a denominator, which improved the computational cost and performance for the KO solution methods during initial testing:
\begin{align}
    r & = \frac{1}{d} \\
    \dot{r} &= - \frac{1}{d^2} \dot{d} \nonumber \\
    &= - r^2 \dot{d} \nonumber \\
    &= r^2 v_A \cos v \label{eq:r_dyns}\\
    \dot{\alpha} &= v_A r \sin v - u \label{eq:alpha_dyns}
\end{align}

Let the state of the system be defined as $x := (r, \alpha)$. The terminal cost and running cost, respectively, are 
\begin{align}
    g \left(x_T\right) &= r_T \cos \alpha_T \\
    h \left(x\right) &= 0.1 r \cos \alpha
\end{align}

The agent wants to minimize the objective while the turret wants to maximize it. This objective represents the fact that the agent wants either to get close to the turret and behind the turret ($\cos \alpha < 0$ and large $r$) or to get as far away as possible if it is in the turret's sights ($\cos \alpha > 0$ and small $r$). The turret, on the other hand, has no direct control of $r$, so it only tries to send $\alpha \rightarrow 0$ (i.e., trying to get the agent in its sights). The running cost captures the fact that both the agent and the turret care about the intermediate steps between the initial and final conditions.

In principle, it might also be possible for the agent to reach the turret ($d = 0$, or $r \rightarrow \infty$). This is a problem both numerically and conceptually, so we also impose a minimum approach distance $d = r = 1$ on the problem. Practically speaking, this minimum approach distance may represent a physical barrier that the agent is unable to cross.

The full specification of the game is therefore
\begin{align}
    J(x, u, v) =  &\ r_T \cos \alpha_T + \int_0^T 0.1 r(t) \cos \alpha(t) \mathop{\mathrm{d} t} \label{eq:objective} \\
    \dot{r} &= r^2 v_A \cos v \\
    \dot{\alpha} &= v_A r \sin v - u \\
    -1 &\leq u \leq 1 \\
    r &\leq 1
\end{align}

The Value function is defined as

\begin{equation}
    \label{eq:Value_definition}
    V(x) = \min_{v(\cdot)} \max_{u(\cdot)} J = \max_{u(\cdot)} \min_{v(\cdot)} J
\end{equation}

\noindent and obeys the saddle-point equilibrium property:

\begin{equation}
    \label{eq:saddle-point-definition}
    J(x, u, v^*) \le \underbrace{J(x, u^*, v^*)}_{V(x)} \le J(x, u^*, v) \qquad \forall u \in U,\ \forall v \in V
\end{equation}

\noindent where $U$ and $V$ are the sets of admissible control strategies for the turret and agent, respectively, and $u^*$ and $v^*$ are their respective equilibrium strategies. For the parameter settings of $v_A = T = 1$, this game was solved in~\cite{vonmoll2024constrained} via the Method of Characteristics (c.f.~\cite{isaacs1965differential,pachter2024synthesis}). The Value function, equilibrium agent strategy, and equilibrium flowfield could not be derived analytically, but these were constructed via backwards integration of the equilibrium state and adjoint dynamics from the terminal manifold.

\subsection{Complementarity-Based Solution with EDMD}
\label{sec:edmd_implementation}




Because we want to know the evolution of the state vectors, we include the identity observable $g_I(x) = x = [r \; \alpha]^\top$ in the dictionary $\Psi$. In addition, we include $g(x) = \cos \alpha$ to approximate the objective function of the turret defense game. For the additional ``lifting'' dictionary functions, we use Random Fourier Features (RFFs)~\cite{rahi07a}
\begin{align}
    \psi_i(x) = \cos(\varphi_i^\top x + b)
\end{align}
\noindent where each unique RFF observable $\psi_i$ is defined by a random vector $\varphi_i$, which we take to be sampled from a diagonal Gaussian distribution (i.e., it has no cross covariance). We found that the performance of the EDMDc model does not depend heavily on the variance chosen for the diagonal Gaussian. For our experiments we used a high variance of $100.0$. The offsets $b$ are uniformly distributed between $0$ and $2\pi$. 

For the turret defense game, note that the original formulation of the state dynamics in Eq.~\ref{eq:r_dyns} and Eq.~\ref{eq:alpha_dyns} have a nonlinear dependence on the controls. As detailed in~\ref{app:edmd}, there are different choices one can make in how to parameterize the control inputs for the EDMDc algorithm to approximate Koopman operators for control systems.

Unfortunately, applying an EDMDc model with an additive dependence on the unlifted controls, as in Section~\ref{sec:edmd_methodology}, to the original formulation of the turret defence game to learn the controlled dynamics did not produce good results -- see~\ref{app:edmd} for details. An equivalent formulation that worked much better broke the agent's control policy into the components of its heading that are parallel and perpendicular to its position relative to the turret:
\begin{align}
    v_A \cos v &\rightarrow v \\
    v_A \sin v & \rightarrow v_{\perp} \\
    v^2 + v_\perp^2 & \leq v_A^2
\end{align}

Technically, we should have $v^2 + v_\perp^2 = v_A^2$. However, relaxing this to an inequality makes the problem easier to solve, and we can show that, for an optimal solution, the inequality will always be at its bound. This changes our dynamics to
\begin{align}
    \dot{r} &= r^2 v \\
    \dot{\alpha} &= r v_\perp - u \\
    v_A^2 &\leq v^2 + v^2_\perp
\end{align}

One additional transformation of the control variables produced even better predictive results:
\begin{align}
    \nu &= r^2 v \\
    \nu_\perp &= r v_\perp \\
    \nu^2 + r^2 \nu^2_{\perp} &\leq v_A^2 r^4
\end{align}

These lifted controls have a one-to-one mapping back to the original agent control variables, so there is no problem in recovering the original control policy. These modifications can be considered as a kind of domain knowledge incorporation -- a lifting of the control inputs that we expect to be advantageous because of our knowledge of the original dynamical system.

To learn the KO representation, we used RFFs trained over the domain 
\begin{align}
    \alpha &\in \left[0,\pi\right] \\
    r &\in \left[0,1\right] \\
    u &\in \left[-1,1\right] \\
    \nu = \cos v, \nu_\perp = \sin v, v &\in \left[0,2\pi\right]
\end{align}

\noindent using gridded data points. To reflect both the symmetry of the problem and the restriction of our domain of training data to $\alpha \geq 0$, we imposed constraints on $\alpha$
\begin{gather}
    0 \leq \alpha_t \leq \pi
\end{gather}
\noindent to respect that limitation. We already have a constraint $r \leq 1$, and the constraint $r \geq 0$ is effectively enforced by the agent's velocity constraints: at $r = 0$, we have
\begin{align}
    0 & \geq \nu_t^2
\end{align}
\noindent so the component of the agent velocity at $r = 0$ perpendicular to that boundary is 0. These constraints on $\alpha$ prevent unnecessary extrapolation of the Koopman dynamics beyond the space of training data (e.g., as might happen during the solution process), and thereby avoid the introduction of additional numerical error on top of normal approximation error.

The restriction of $\alpha$ to be non-negative allowed us to utilize the symmetry of the underlying problem, as has been done previously in the Koopman literature~\cite{sinha2020koopman}. More specifically, for any trajectory $\left\{\alpha,r,u,\nu,\nu_{\perp}\right\}_t$, there is an equivalent trajectory $\left\{-\alpha,r,u,\nu,-\nu_{\perp}\right\}_t$. This means that we only need to learn a KO model for $\alpha \geq 0$ to cover $\alpha \in \left[-\pi,\pi\right]$. We also have the equivalency $\left\{\alpha + 2n \pi,r,u,\nu,\nu_{\perp}\right\}_t, n \in \mathbb{Z}$, but this is trivial for the problem under consideration. The reduction of the problem to non-negative $\alpha$, on the other hand, allows us to sample more densely from the domain and thereby, in principle, to get a better EDMD representation without additional cost.

As an additional bit of domain knowledge incorporation, instead of learning a general data-driven form for $Q$, we explicitly specified $\cos \alpha$ as an observable, as mentioned previously, which allowed us to specify the form of 
\begin{gather}
     Q \equiv Q_g = \frac{1}{0.1 \Delta t} Q_h 
\end{gather}
\noindent \textit{a priori}: a value of 1.0 in the entry of $Q$ corresponding to $r \cos \alpha$ and zeros elsewhere.

Note that the turret and agent optimizations are slightly different because of the constraints that apply to the agent but not the turret -- in particular, the constraint on $r$ is not relevant to the turret because the turret has no control over $r$. This is particularly important when we look at a subset of the KKT conditions for the turret
\begin{align}
    0 =& -0.1 \Delta t \left[\left(Q_{\psi,\psi} + Q^T_{\psi,\psi}\right) \psi_t + \left(Q_{\psi,1} + Q^T_{1,\psi}\right)\right] \nonumber \\
    &+ \left(\sigma_{\pi,t} - \sigma_{0,t} \right) \circ c_{\alpha,select} - \lambda_{u,t-1} + \lambda_{u,t}^T A  & t = 2,\ldots,T-1 \\
    0 =& - \left[\left(Q_{\psi,\psi} + Q^T_{\psi,\psi}\right) \psi_T + \left(Q_{\psi,1} + Q^T_{1,\psi}\right)\right] \nonumber \\
    &+ \left(\sigma_{\pi,t} - \sigma_{0,t} \right) \circ c_{\alpha,select} - \lambda_{u,T-1} & t = T
\end{align}
\noindent and the corresponding subset of the KKT conditions for the agent
\begin{align}
    0 =& 0.1 \Delta t \left[\left(Q_{\psi,\psi} + Q^T_{\psi,\psi}\right) \psi_t + \left(Q_{\psi,1} + Q^T_{1,\psi}\right)\right] \nonumber \\
    &+ \left(\sigma_{\pi,t} - \sigma_{0,t} \right) \circ c_{\alpha,select} \nonumber \\
    & + \left(\sigma_{r,t} + r_t \nu^2_{\perp,t} \mu_t - 2 v_A^2 r_t^3 \mu_t\right) \circ c_{r,select} \nonumber \\
    &- \lambda_{v,t-1} + \lambda_{v,t}^T A  & t = 2,\ldots,T-1 \\
    0 =& \left[\left(Q_{\psi,\psi} + Q^T_{\psi,\psi}\right) \psi_T + \left(Q_{\psi,1} + Q^T_{1,\psi}\right)\right] \nonumber \\
    &+ \left(\sigma_{\pi,T} - \sigma_{0,T} \right) \circ c_{\alpha,select} \nonumber \\
    & + \left(\sigma_{r,T} + r_T \nu^2_{\perp,T} \mu_t - 2 v_A^2 r_T^3 \mu_T\right) \circ c_{r,select} \nonumber \\
    &- \lambda_{v,T-1} & t = T \\
    0 \leq& 1 - r_t \perp \sigma_{r,t} \geq 0 & t = 2,\ldots,T \\
    0 \leq& v_A^2 r^2 - \nu^2 - r^2 \nu^2_\perp \perp \mu_t \geq 0 & t = 1,\ldots,T-1
\end{align}
\noindent where $c_{r,select}$ is a vector of zeros with a 1 in the entry corresponding to $r$, $c_{\alpha,select}$ is a vector of zeros with a 1 in the entry corresponding to $\alpha$ and
\begin{align}
    0 &\leq \pi - \alpha_t \perp \sigma_{\pi,t} \geq 0 & t = 2,\ldots,T &  \\
    0 &\leq \alpha_t \perp \sigma_{0,t} \geq 0 & t = 2,\ldots,T &
\end{align}
\noindent are the complementarity conditions corresponding to the bounds on $\alpha$. Note that in this case, $\lambda_u$ is more than just $-\lambda_v$ because the agent has constraints that the turret does not. Lifting the agent's control variables improves the EDMD representation, but as the agent's KKT conditions show, it makes the MCP nonlinear, and thus harder to solve. As a nonlinear MCP, the initial conditions provided to the algorithm matter, and thus we provide initial guesses based on simple heuristics---the turret turns at full speed towards the agent while the agent moves directly toward the turret for $\alpha \geq 1.0$ or perpendicular to the turret for $\alpha < 1.0$. Regardless, we solve this MCP using the PATH algorithm~\cite{dirkse1995path} and the MCP capability of JuMP~\cite{Lubin2023}.

\subsection{Resolvent-Based Solution}
\label{sec:resolv_implementation}

\subsubsection{Finite Basis Selection}

For the resolvent-based approach, we used the same domain as specified for the original problem and MCP formulation. Initially, we attempted to use ApproxFun.jl~\cite{ApproxFun.jl-2014} as an adaptive way to control the resolvent calculation error throughout the optimization process. This required calculating derivatives through the resolvent and explicitly handling constraints (see~\ref{app:resolve}). We were able to do both of these things, but the optimization process was very slow and struggled to converge. As such, we instead used a grid of Radial Basis Functions (RBFs) over the domain to give us a fixed-size basis. We actually spread the RBF nodes slightly outside of the domain to get better function approximation within the domain and then selected a denser grid of evaluation points within the domain. This allowed us to precalculate the necessary kernels once and then use them repeatedly in the optimization. We used a 5x5 grid of RBF kernels on $\left(r,\alpha\right) \in \left[-0.1,1.1\right] \times \left[-0.1 \pi, 1.1 \pi \right]$ with a 25x25 grid of evaluation points on $\left(r,\alpha\right) \in \left[0,1\right] \times \left[0, \pi \right]$.

\subsubsection{Constraint Enforcement}

When implementing the resolvent-based approach, we found that splitting the agent's control policy into $v$ and $v_{\perp}$ components (as discussed in Section~\ref{sec:edmd_implementation}) significantly improved convergence of the optimization process. This meant, however, needing to impose the
\begin{gather}
    v^2 + v_\perp^2 \leq v_A^2
\end{gather}
\noindent constraint. It is relatively straightforward to impose a constraint like this directly on control inputs, at least approximately, by imposing those constraints at each evaluation point within the domain (as described in Section~\ref{sec:resolve_methodology}); we could not impose this constraint as an equality because it would overconstrain the problem. However, we found that, due to the limits of the finite-basis approximation, $v^2 + v_{\perp}^2$ would sometimes be small relative to $v_A^2$ -- especially when $r$ was small. To improve solution quality, therefore, we also imposed a lower bound
\begin{gather}
    v^2 + v_\perp^2 \geq 0.9 v_A^2
\end{gather}
\noindent on the system, but using $\nu = r^2 v$ and $\nu_\perp = r v_\perp$ was not necessary.

Imposing constraints on the states (e.g., $r \leq 1$) was more complicated -- especially when the constraint is supposed to hold at each point in time and the agent is incentivized, by the objective, to violate it. In principle, it would be possible to do this by using the resolvent to calculate $r_t$ at a sufficiently dense set of intermediate time points and imposing a constraint on the calculated $r_t$ values from each evaluation point in the domain. However, this would greatly add to the size and complexity of the optimization problem, making it much more computationally expensive to solve. Even just imposing this kind of constraint on a terminal condition $r_T \leq 1$ could be problematic. 

However, in this particular case, these was a far simpler alternative approach. If we imposed the boundary condition $\dot{r} \leq 0$ on the boundary $r = 1$, this ensured that the agent never crossed the boundary. Given our definition above, this amounted to specifying $v \leq 0$ on $r = 1$. To enforce this, we simply defined a set of evaluation points restricted to $r = 1$ and proceeded as before. We could do this because the agent can always choose a control policy that satisfies the constraint -- it does not create any feasibility problems for the optimization. If instead we had something like
\begin{align}
    \dot{r} & = r^2 v + \alpha
\end{align}
\noindent this would not necessarily work -- if $\alpha = \pi$ and $v_A = 1$, for example, then $\dot{r} \geq 0$ for any value of $v$ that satisfies $v^2 + v_\perp^2 \leq v_A^2$. The domain boundaries in this problem did not present any issues, though -- partially because of how we chose the value of $v_A$ and the bounds on the turret's maximum angular speed. To match the EDMD-MCP formulation, and to similarly avoid numerical problems with extrapolation, we enforced the constraints $\dot{\alpha} \leq 0$ at $\alpha = \pi$ and $\dot{\alpha} \geq 0$ at $\alpha = 0$; again the existing constraint on the agent velocities effectively enforced $\dot{r} \geq 0$ at $r = 0$ as noted in Section~\ref{sec:edmd_implementation}. Unlike the constraint at $r = 1$, both the agent and the turret could influence $\dot{\alpha}$ directly. Therefore, we enforced the constraints
\begin{align}
    0 &\geq v_A r v - u  \\
    0 & \leq v_A r v - u
\end{align}
at $\alpha = 0$ and $\alpha = \pi$, respectively, using an approach analogous to the one used on $\dot{r}$ at $r = 1$.

\subsubsection{Numerical Optimization Procedure}

With the problem formulated in Nonconvex~\cite{MohamedTarekNonconvexjl}, we used IPOPT~\cite{wachter2006implementation} to alternate between solving the agent and the turret optimization problems using warm starts from each previous iteration. We used
\begin{align}
    u \left(r,\alpha\right) &= 1 - \frac{\left(\alpha - \pi\right)^4}{\pi^4} \\
    v \left(r,\alpha\right) &= \cos \left(\pi - 0.5 r \alpha\right) \\
    v_{\perp} \left(r,\alpha\right) &= \sin \left(\pi - 0.5 r \alpha \right)
\end{align}
\noindent as initial solutions from which to work. These solutions are reasonable given the goals and constraints of the players, but they are not optimal. The alternating optimizations were only allowed to run for 100 iterations before switching to the other player. We found that limiting the number of iterations and performing multiple sequential optimizations (using warm starts with small perturbations) significantly helped with convergence -- especially for the agent, who had a more nonlinear problem to solve.

Further details of the numerical implementation for the resolvent-based method are provided in~\ref{app:resolve}.

\section{Results}
\label{sec:results}

\subsection{Analytical Solution}

As stated previously, the analytic solution presented in~\cite{vonmoll2024constrained} corresponds to the parameter settings $v_A = T = 1$; the results are briefly summarized here.
Note that the case when $v_A \neq 1$ significantly complicates the analysis of the trajectories for which the constraint ($r = 1$) becomes active and thus was not considered in the prior work. Within the analytic solution, there are two singular surfaces (c.f., e.g.,~\cite{basar1982chapter}) namely the Universal Line (UL) and the Dispersal Line (DL) given, respectively, by
\begin{align}
    \label{eq:UL_and_DL}
    \mathcal{U} &:= \left\{ (r,\alpha) \mid \alpha = 0,\ r \in \left( 0, 1 \right] \right\} \\
    \mathcal{D} &:= \left\{ (r, \alpha) \mid \alpha = \pi,\ r \in \left( 0, 1 \right] \right\}
\end{align}
These lines are singular because the equilibrium control strategies of the players cannot be uniquely determined from the necessary conditions for optimality. In both cases, symmetry is to blame (i.e., since $\alpha \in \mathbb{S}^1$). On the UL, $\alpha = 0$ and thus the turret's look angle is directly aligned with the agent's position.
With $r \in (0, 1]$, it is clear that the turret has superior control authority over the $a$ state and thus the agent has no hope of improving its outcome by deviating from $\alpha = 0$. Thus the equilibrium turret and agent strategies on the UL are to not rotate ($u = 0$) and run directly away ($v = \pi$), respectively. On the DL, $\alpha = \pi$ and thus the turret can either begin turning clockwise \textit{or} counterclockwise in order to turn ``toward'' the agent. Depending on the choice of the turret, the agent has a corresponding equilibrium heading that involves some component of velocity that heads in the same direction the turret turns. These singular lines cause the solution to be symmetric about $\alpha = 0$ and $\alpha = \pi$, and thus, without loss of generality, one need only consider $\alpha \in \left[0, \pi \right)$.

The overall structure of the solution can be summarized by describing the three fundamentally different types of equilibrium trajectories.
ULs, by their nature, are, themselves, an equilibrium trajectory.
Additionally, some trajectories beginning elsewhere in the state space feed into the UL. 
\begin{definition}[Universal Line Tributaries]
   Equilibrium trajectories that begin inside the state space (i.e., $r \in \left( 0,1 \right)$ and $\alpha \in \left( 0, \pi \right)$) and reach the UL, with $\alpha = 0$ within the time horizon of the game.
\end{definition}
For the special case of $v_A = 1$ the analysis of the constraint is as follows. When $r = 1$, the agent is at the closest allowable distance to $T$ and their angular speeds are equal. Thus the agent can ensure $\dot{\alpha} = 0$. When $\alpha > \pi/2$, the running cost, $r \cos \alpha < 0$ and thus the agent, being the minimizer has incentive to maintain this state (i.e., $\dot{r} = \dot{\alpha} = 0$). Alternatively, when $\alpha < \pi/2$, there is no incentive for the agent to remain on the constraint. Thus any trajectory for which the constraint becomes active remains on the constraint.
\begin{definition}[Constrained Trajectories (for $v_A = 1$)]
    Equilibrium trajectories that may begin away from the constraint ($r < 1$), reach the constraint ($r = 1$) within the time horizon of the game, and subsequently remain there.
\end{definition}
However, it may be the case that the state remains in the ``interior'' of the state space. That is, the state neither reaches the UL nor the constraint boundary.
\begin{definition}[Regular Trajectories]
    Equilibrium trajectories for which $\alpha(t) > 0$ for all $t \in [0, T]$ and $r(t) \neq 1$ for all $t \in [0, T]$.
\end{definition}

\subsection{Computational Solutions}
\label{sec:comp res}


For the given differential game, we only needed to train a single EDMDc model. This could involve hyperparameter optimization, such as the size of the RFF dictionary or the distribution the RFF frequencies are sampled from. However, because EDMDc models are learned via least-squares regression, they were efficient to train in this case -- at most a few minutes on a laptop. We tested the predictive accuracy of a trained EDMDc model by feeding it sample control inputs and computing the error between the EDMDc state trajectories and those from numerically integrating the sample controls in the original equations for the state dynamics. 

With an EDMDc model in hand, equilibrium solutions for the game were approximated by solving the MCP \emph{for each individual initial condition} -- in contrast to the resolvent-based method that globally solves for all equilibrium solutions. For the turret defense game, equilibrium trajectories were computed individually for given initial states $(r_0, \alpha_0)$ and time horizon $T$. In order to help the MCP solver converge, we provided initial solution guesses using simple heuristics based on system symmetries and instantaneous objectives (e.g. if the agent should turn toward or away from the turret at each time), as discussed in Section~\ref{sec:edmd_implementation}.

We used the EDMDc-based MCP method to approximate equilibrium trajectories for a little over $1500$ initial conditions using $T = 1.0$. There were around $900$ Regular solutions, and the rest split evenly between Constrained and Universal Line. In all cases, the MCP solver converged to a local optimum using our initial guess heuristics. The time to solution varied for each initial condition, but it was generally very fast. The median solution time was just under $20$ seconds, and there were seven outliers that took over an hour, with the longest just over three hours. 

The computational costs associated with the resolvent-based solutions were, in many ways, the exact opposite of the EDMDc-based MCP method: there was a single, expensive set of optimizations to calculate the solutions for the entire domain, and evaluating the solutions then just involved integrating an ODE instead of solving an MCP. We let the optimizations alternate back and forth 100 times, and by the end, each optimization would only change the objective function by about 0.1\% over the course of a few iterations before reaching optimality; the process essentially converged to an equilibrium. It took almost seven hours to do so, but once it was done, any given solution in the domain could be calculated quickly and easily.

\begin{figure}[ht]
\centering
\includegraphics[width=0.48\textwidth]{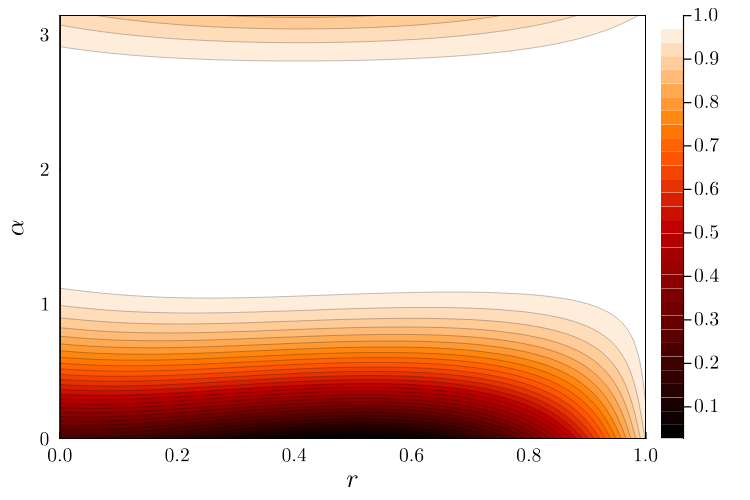}
\includegraphics[width=0.48\textwidth]{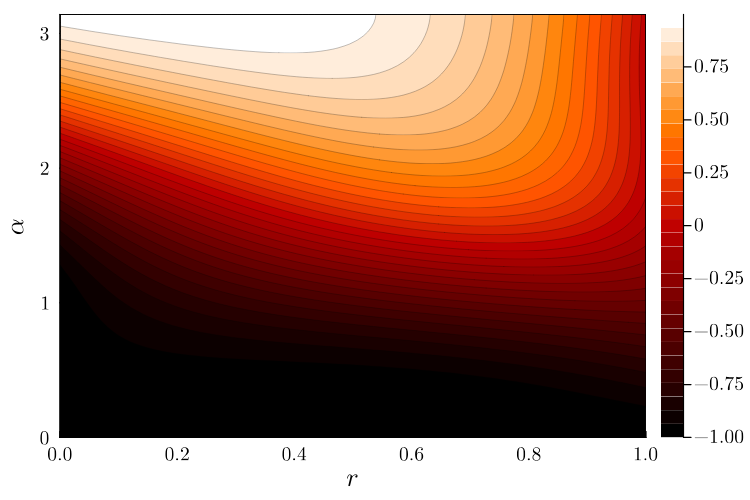}
\includegraphics[width=0.48\textwidth]{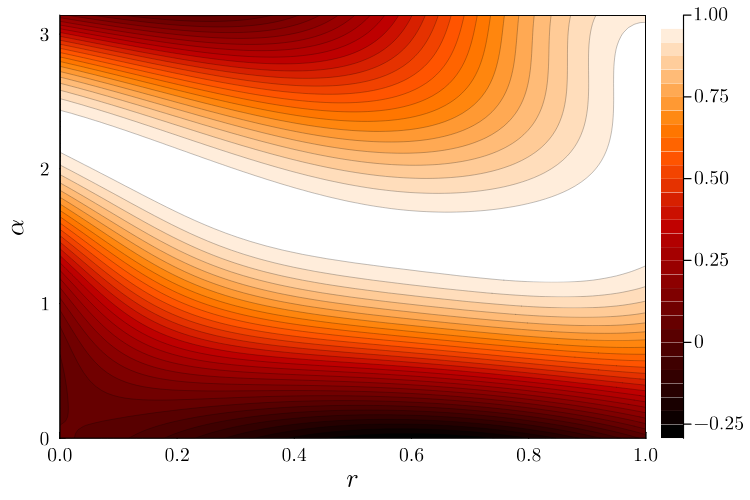}
\includegraphics[width=0.48\textwidth]{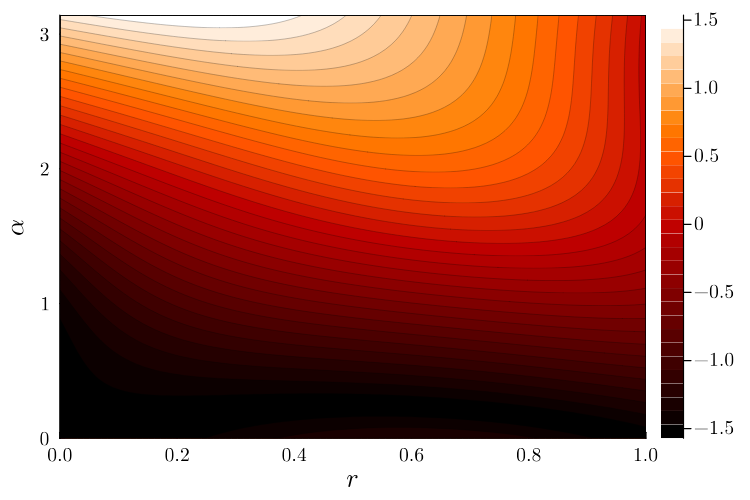}
\caption{Resolvent-Based Equilibrium Control Policy Solutions: $u$ (Top Left), $v$ (Top Right), $v_\perp$ (Bottom Left), Agent's Relative Heading (Defined as $v$ in Figure~\ref{fig:turret-defense-schematic}, Bottom Right)}
\label{fig:resolv_policies}
\end{figure}

The feedback control policies calculated by the resolvent-based approach are shown in Figure~\ref{fig:resolv_policies}. The turret's control policy $u$ approximates what we know the optimal policy to be (namely, $u = 1$ everywhere with a sharp discontinuity at $\alpha = 0$~\cite{vonmoll2024constrained}). Over most of the domain, the approximation was fairly good. However, the fixed-basis approach had less fidelity around the $\alpha = 0,\pi$ boundaries. In particular, the feedback policy was necessarily smooth (as the weighted sum of smooth functions), and its ability to produce a sharp drop at $\alpha = 0$ was limited.

The agent's policy also qualitatively matched the known (non-feedback) equilibrium solution but struggled around discontinuities. For the agent, the key discontinuity is around $\alpha = 1$ for large values of $r$. In that region, the agent would head towards the turret (i.e., $v > 0$). However, to avoid violating the $r \leq 1$ constraint, there is a boundary such that, for some $\epsilon$, $v_{\perp} \geq \epsilon > 0$ on one side of the boundary and $v_{\perp} \leq - \epsilon < 0$ on the other side. The bifurcation in the agent's control policy is difficult to capture with the basis functions used here. This highlights one of the limitations of a feedback control policy. For the analytical and EDMD-MCP solutions, the agent's behavior in that region depends on the time at which the agent enters that region. A feedback control policy cannot take that into consideration, however, as it has no explicit time dependence.

It is also worth noting that for both the agent and the turret, there is little incentive to produce highly accurate solutions for small values of $r$. Since the objective is essentially proportional to $r$, the contributions of small-$r$ points to the overall objective $J$ are small, and thus deviations from optimality in those regions will, in general, have small impacts on the overall payoff.

\begin{figure}[ht]
\centering
\includegraphics[width=.48\textwidth]{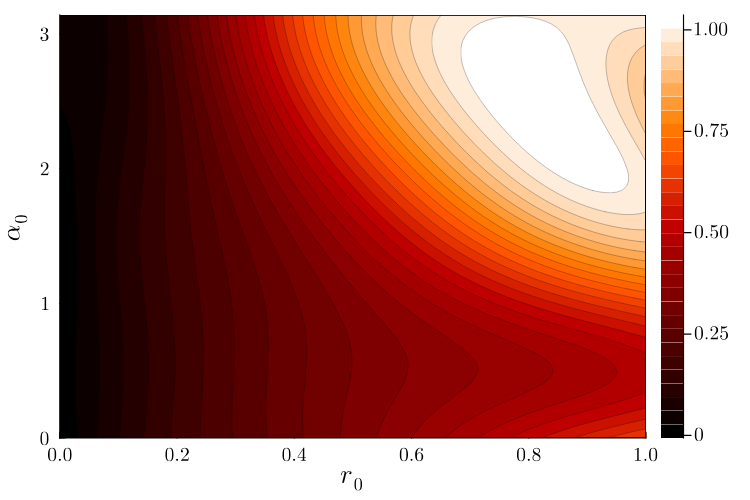}
\includegraphics[width=.48\textwidth]{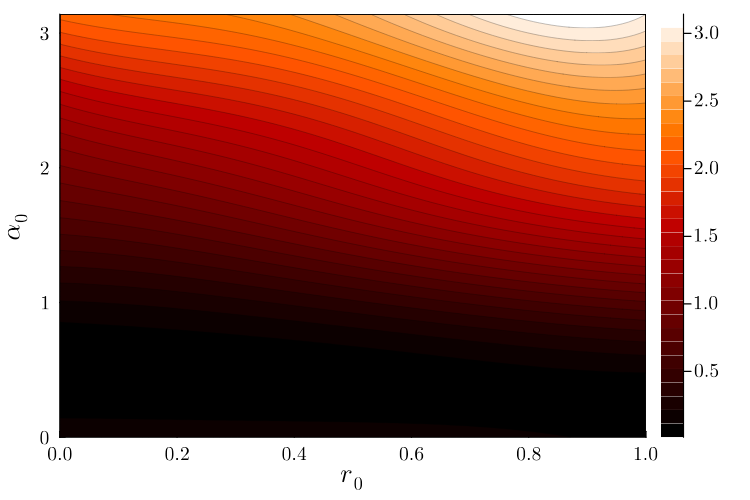}
\caption{Resolvent-Integrated Results: Final $r$ Value (Left), Final $\alpha$ Value (Right)}
\label{fig:final results}
\end{figure}

Figure~\ref{fig:final results} shows the final $r$ and $\alpha$ values given the calculated resolvent-based control policies. Note that here, the axes are $r_0$ and $\alpha_0$, indicating initial conditions, whereas the control policies are purely functions of the state (not the initial state), so that axes are $r$ and $\alpha$. In principle, for the trajectories starting in the upper righthand part of the domain, we would expect $r_T$ to be exactly 1: once the agent hits $r = 1$, it can maintain that position and is incentivized to do so as long as $\alpha > \pi/2$. However, the previously noted trajectory bifurcation for the agent combined with the inherent limitations of this particular finite basis and a feedback control policy in general mean that there is a small blip around $\alpha = 3 \pi/4$ where the agent leaves the $r = 1$ surface. The $\alpha$ trajectories are a bit better behaved in this regard: as we would expect, $0 \leq \alpha_T \leq \alpha_0$.

\begin{figure}
\centering
\includegraphics[width=0.62\textwidth]{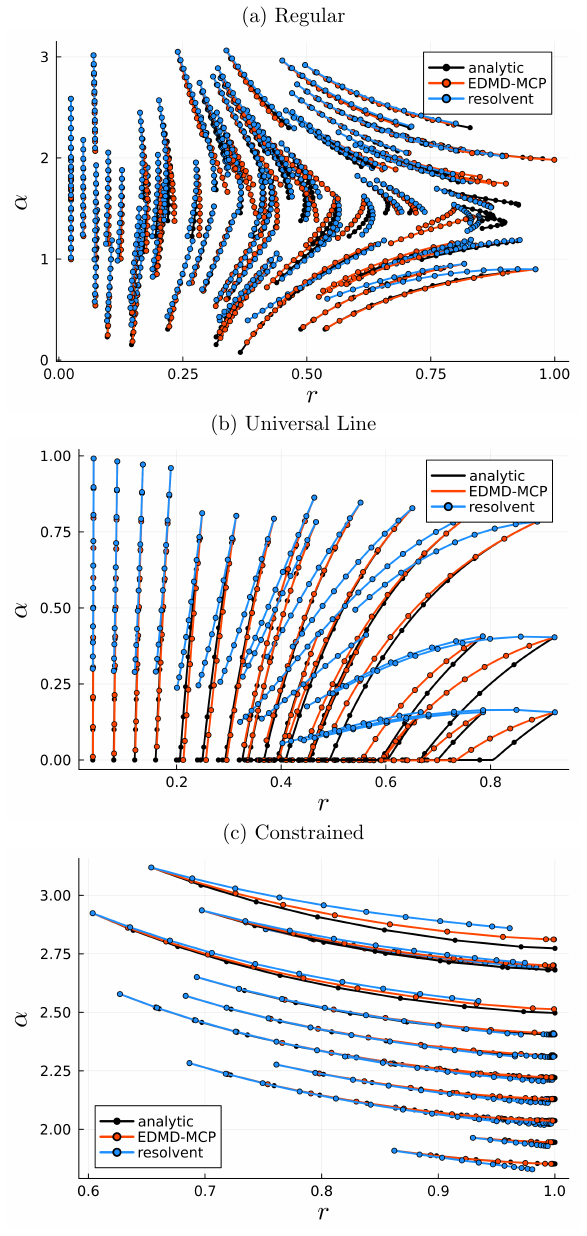}
\caption{Equilibrium Trajectory Comparisons}
\label{fig:trajectories}
\end{figure}

Figure~\ref{fig:trajectories} compares EDMD-MCP, resolvent, and analytic equilibrium trajectories for a sample of initial conditions. For Regular (a) and Constrained (c) trajectories, there was general qualitative agreement in trajectories. For most Regular and Constrained trajectories, there was strong quantitative agreement, too, with the trajectories largely overlapping. We see the most error with Regular solutions for $r > 0.75$ and $\alpha$ around $\pi/2$ -- right where the agent tends to shift behavior from approaching the turret to running away from the turret. The agent tends to change direction and run away (decreasing $r$) too quickly with our numerical approximations, particularly the EDMD-MCP. Generally, the EDMD-MCP solutions are closer to the analytical solutions than the resolvent solutions are, and this holds true particularly for $\alpha < 1$.

The Universal Line trajectories were the most difficult to approximate. The EDMD approximation of the dynamics of the $\cos \alpha$ observable used for the objective function did not fully capture the nonlinearity of this observable, which penalizes the agent for not running directly away when the turret is pointed at it directly. The explicit constraints $0 \leq \alpha \leq \pi$ in the MCP helped to alleviate this problem. Similar to the explicit constraints on $r$, this forced the agent to remain within the set boundaries of the problem. However, these boundaries on $\alpha$ are based on rotational symmetries of the problem, and so violating the boundary constraints would still give physical solutions, even if they were not the correct equilibrium solutions. By contrast, we have found that the MCP solutions, when using a bad EDMDc model, could produce unphysical trajectories with $r < 0$ without using explicit constraints on $r$. 
Therefore, although boundary constraints in the MCP may not be required in principle, given a perfect EDMDc model, in practice we found them necessary to ensure good quality equilibrium solutions. With those constraints in place, EDMD-MCP matched the Universal Line analytical solutions qualitatively over the whole domain, with better agreement for smaller $r$ values. The resolvent solutions did not match the analytical solutions nearly as well because of the inability to produce $u=1$ solutions near $\alpha = 0$, as shown in Figure~\ref{fig:resolv_policies}.

\begin{figure}
\centering
\includegraphics[width=0.62\textwidth]{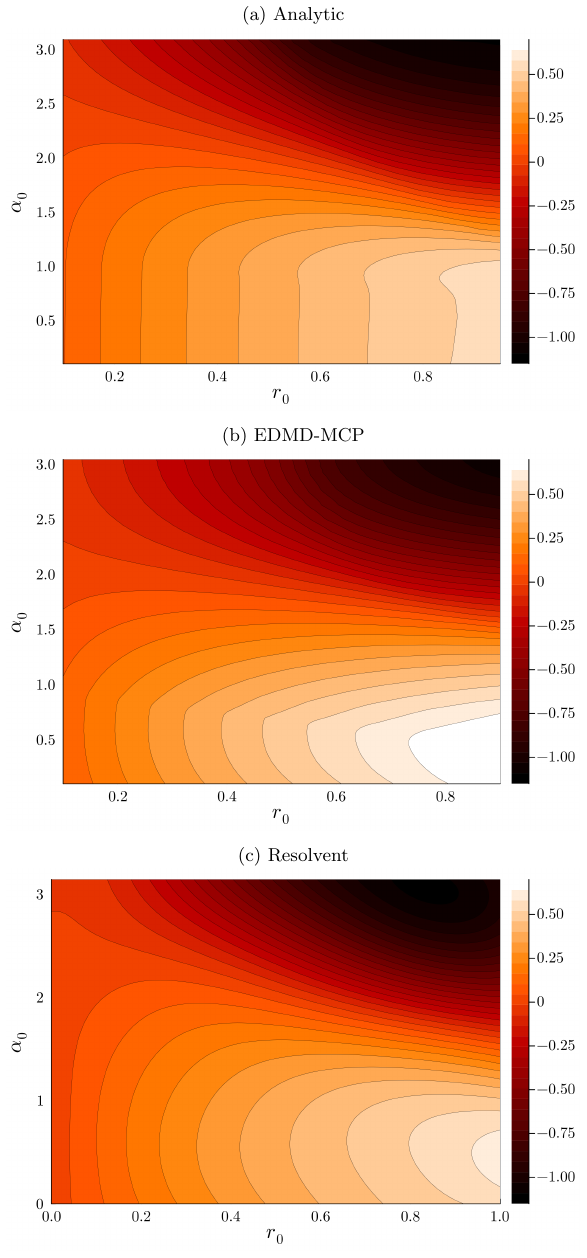}
\caption{Values of the Turret Defense Game}
\label{fig:game_values}
\end{figure}

Finally, Figure~\ref{fig:game_values} compares the value of the game as a function of the initial conditions $r_0$ and $\alpha_0$ for the analytic, resolvent, and EDMD-MCP solutions. The EDMD-MCP and analytic game values are interpolated over the full $r$-$\alpha$ domain using the roughly $1500$ sample trajectories. Qualitatively, we see strong agreement between all three, although some details differ. In these plots, lighter regions favor the turret, while darker regions favor the agent. The most notable feature is that the EDMD-MCP solutions overly favor the turret in the bottom right region on the phase space, whereas the resolvent favors the agent slightly for small $\alpha$ values -- again, because of the basis limitations in generating the feedback control policy.

\section{Discussion}
\label{sec:discussion}

\subsection{Comparing and Contrasting Complementarity- and Resolvent-Based Methods}

Complementarity- and resolvent-based solution methods have very different properties and thus would likely be used for different things. The resolvent-based approach is limited to a feedback control policy, which is a restriction on the set of available control policies, but it also allows us to solve for the entire domain simultaneously rather than having to solve the system one trajectory at a time. The computational cost of solving for an equilibrium control policy over the entire domain is significant, but once that is done, querying the saved control policy has a negligible computational cost. 

The assumption of a feedback control policy also allows us to avoid the linear control assumption used in the MCP -- i.e., that we can separate the lifted state and control terms out so that the lifted state dynamics are linear in the (potentially unlifted) control terms. Avoiding the need for such an assumption suggests that, compared to the complementarity-based approach, the resolvent-based approach will have a higher level of fidelity in the solutions it calculates. The tradeoff is that the resolvent does not gain the computational benefits that the linear control assumption provides: we cannot use KKT conditions to solve directly for the equilibrium, and the optimizations that we perform are necessarily nonlinear. Furthermore, as we saw in the results presented above, the basis can strongly limit the global approximation of an equilibrium feedback policy.

Alternatively, the feedback control policy of the resolvent-based approach means that it may not be necessary to lift the control terms to get good results. In the turret defense problem considered here, for example, lifting the control inputs using prior knowledge of the system proved to be necessary for the discrete-time EDMD approach to approximate the dynamics well. That was not a requirement for the complementarity-based approach, strictly speaking, but rather a result of trying to learn the KO representation of a dynamical system with a highly nonlinear control function; under other circumstances, such efforts might not be necessary. Not needing to lift the control inputs reduces the effort required to implement the resolvent-based approach.

In its current implementation, the resolvent-based approach only works on a continuous-time, white-box representation of the system in question. It might be possible to adapt the method to deal with a black-box, discrete-time system -- many real-world systems do not have known ODE representations and instead only produce sampled data -- but this would require an extra layer of computation or approximation to turn that data into a continuous-time form. Finally, because we are solving over the entire domain, imposing constraints on the states becomes more complicated. In this turret defense case, we were able to specify state constraints implicitly (i.e., $r \leq 0$) by imposing constraints on the control inputs on the boundary of interest. That may not always be possible, though, and in that case, it may be necessary to relax the constraints and/or to define the domain of interest in a highly selective fashion. We are also not directly solving for the Nash equilibrium. Rather, we are using alternating optimizations to converge (approximately) to that point. This worked fairly well in the case presented here, but such approaches are known to suffer convergence issues in other contexts~\cite{schaefer_competitive_2019}.

The MCP-based approach, on the other hand, solves for open-loop control trajectories one trajectory at a time. There is an upfront training process involved, but that training is relatively inexpensive, and solving for a single trajectory also quick. It works in discrete time, can directly apply constraints to states, and may use sampled black-box data. Although it does not, in principle, require knowledge about the system dynamics beyond what the data provides, it is possible to incorporate domain knowledge into the KO representation as we showed here. The linear control assumption may negatively impact the accuracy of the KO representation, but that assumption facilitates the direct solution of the Nash equilibrium rather than iterating back and forth as with the resolvent.

Both approaches are similar in that they use fixed basis approximations, but otherwise, they are essentially mirror opposites. For situations like this, where we are interested in solutions over the entire domain, the resolvent-based approach is probably more useful. As fast as the MCP is, its time-to-solution is still comparable to the duration of the game itself. Once solved, however, the resolvent-based approach would only require quick function evaluations to determine the equilibrium policies for both players. More generally, this characteristic could be advantageous in low-latency systems where there would not be time to solve even a fast MCP problem. As the dimensionality of the problem increases, the use of dimension reduction via, e.g.,  Deep-DMD would give the complementarity-based approach an advantage. The resolvent-based approach requires a basis in which to perform the calculations. If that basis were changing at each optimization iteration, those iterations would then also be accompanied by a large number of overhead calculations that were only performed once in the fixed-basis case.



\subsection{Future Work}

Our results here suggest several different future avenues for research. Firstly, the turret problem described here had some challenging properties (e.g., the bifurcation in feedback policies as discussed in Section~\ref{sec:comp res}) but was still quite small. A key next step, therefore, would be to apply these methods to larger problems. In these larger problems, the complementarity-based approach is likely to be more feasible than the resolvent-based approach, though computational costs will increase for both.  The curse of dimensionality is the main issue for the resolvent-based approach: the size of the basis and the number of sample points needed to represent the entire state space will grow exponentially in the dimension of the problem.  By contrast, Deep-DMD approaches, including Deep-DMD with autoencoders, can handle higher-dimensional dynamics much more flexibly and easily \cite{bakker2020koopman}.  That being said, the computational cost of solving MCPs grows superlinearly in the size of the problem, and if the goal is to solve these kinds of dynamic games in real-time or near-real-time, compact discrete-time KO models will be necessary.  In other words, computational costs will be a major consideration for the KO-based solution of larger games, but the complementarity-based approach may be less vulnerable to the curse of dimensionality.

Secondly, many games do not have fixed time horizons but instead terminate once certain conditions are met. There are ways in which this could be implemented for the EDMD-MCP method without fundamentally changing the approach, but trying to apply this to the calculation of a global feedback solution, where the time horizon varies over the domain, would be somewhat more complicated. See our derivations and sample problem in~\ref{app:resolve} for further consideration of this subject. On a related note, however, extending the resolvent-based methodology to a time-dependent feedback control policy would be relatively straightforward by simply appending the time as an additional state with a constant time derivative of one (e.g., as is done to transform non-autonomous systems into autonomous ones).  Increasing the dimension of the state space will increase the computational cost of solving the problem, as discussed in the previous paragraph, but it would not fundamentally change any of the methodology presented above.

Thirdly, this methodology was developed for two-player zero-sum games. The EDMD-MCP approach is agnostic with regards to the number of players or the relationship between their objectives. In principle, there would be no problem in applying this method to $n$-player, general-sum games. The resolvent-based approach, however would need further modification. In particular, the alternating optimizations approach would likely struggle more to converge in an $n$-player, general-sum context. This does not rule out the use of the resolvent for solving such problems -- technically, the resolvent is only being used to calculate the cost functional -- but other methods for converging to an equilibrium, given the resolvent calculations, may be necessary. 

Because the KG is a function of the control policies, and because the KG shows up within a number of matrix inversions, as part of the resolvent approximation, trying to apply a direct equilibrium solution approach, like the MCP, would be difficult. Calculating the necessary derivatives through all of those matrix inversions would be computationally expensive, and because each player's optimization would likely not be convex, the KKT conditions in which those derivatives show up would not necessarily define a unique equilibrium.

Fourthly, the formulation presented in Section \ref{sec:methodology} assumed that cost was solely a function of the state. This is not necessary in order to employ the proposed methods, but it does raise a greater number of implementation questions.  For the complementarity-based approach, for example, it may be possible to represent control costs via the addition of quadratic terms to the objective:

\begin{gather}
    \max_{u_t} \min_{v_t} \zeta_T^T Q_g \zeta_T + \Delta t \sum \limits_{t=1}^T \left(\zeta_t^T Q_h \zeta_t^T + u_t^T Q_u u_t + v_t^T Q_v u_t\right)
\end{gather}

\noindent where $Q_u$ is negative semidefinite and $Q_v$ is positive semidefinite. This specific implementation will not be applicable in all cases, but given that the norm of the control input is a common cost component in optimal control problems, it will likely be applicable in many cases; the updated optimality conditions are then straightforward to derive.  More complicated approaches that involve lifted control variables in the objective or cross-terms (e.g., $\zeta_t^T Q_{\zeta,u} u_t$) would also be possible. The resolvent-based approach is somewhat simpler: if the running cost is $h\left(x_t,u_t,v_t\right)$ and the feedback control policies are linear functions of the pre-defined bases, then

\begin{align}
    \tilde{h} \left(x,u_k,v_k\right) &= h\left(x,\sum_k u_k \phi_k\left(x\right),\sum_k v_k \phi_k\left(x\right)\right) \\
    J \left(T\right) & \approx \sum_{i,j} \beta_j \left(u_k,v_k\right) \left(I z_j - L\left(u_k,v_k\right)\right)^{-1} \left( g \left(x_i\right) + \frac{1}{z_j}\left[ 1 - e^{-z_jT}\right] h \left(x_i,u_k,v_k\right) \right)
\end{align}

The extension of the proposed method to problems where the objective includes control inputs is therefore not a large one.

Finally, for some classes of differential games, it may be feasible to use EDMD to leverage existing LQ game solution methods instead of using MCPs. Such an approach would naturally lend itself to differential games that lack explicit constraints on system states or player controls; it was not suitable for the turret defense problem considered here because of the constraints inherent to that problem. In some games with such constraints, though, it may be possible to relax them by removing them and replacing them with appropriate penalty functions in the objective. The penalty functions would then need to be quadratic functions of the lifted states. An EDMD-LQ method would have a much lower run-time cost than the EDMD-MCP method, and under some circumstances, it could even provide a global feedback solution (like the resolvent-based approach, but in discrete time). Such an approach might be key to addressing the computational cost challenges identified previously in this section. It would be worthwhile to explore the pros and cons of EDMD-LQ in comparison with EDMD-MCP and the resolvent-based approach.

\section{Conclusions}
\label{sec:conclusions}

In this paper, we defined, tested, and discussed two fundamentally different ways of using the Koopman operator to solve two-player zero-sum differential games. The first way uses the resolvent to produce a continuous-time feedback solution over the entire problem domain via the use of alternating optimizations. The second way uses a discrete-time Koopman approximation and solves directly for the Nash equilibrium on a trajectory-by-trajectory basis using optimality conditions. Both methods were generally successful at doing so for the turret defense problem presented here, but the process highlighted the relative advantages and disadvantages of each method. Future work includes applying these methods to larger and more general differential games.

\backmatter

\bmhead{Acknowledgment}
The authors would like to thank Andrew Horning for his insightful discussions on resolvent-based contour methods for approximating Koopman semigroups. Additionally, The authors would like to acknowledge the support of Air Force Office of Scientific Research Lab Tasks \#21RQCOR083 and \#24RQCOR001. These tasks are jointly funded by Drs.~Fariba Fahroo, Frederick Leve, and Warren Adams.  DISTRIBUTION STATEMENT A. Approved for public release. Distribution is unlimited. AFRL-2024-5791; Cleared 15 Oct 2024.

\section*{Declarations}

The authors have no competing interests to declare that are relevant to the content of this article.



\appendix
\section{Extended Dynamic Mode Decomposition}
\label{app:edmd}

Consider first an autonomous dynamical system without control inputs, $\dot{x} = f(x)$. The dynamics $f$ define a semigroup of flow maps $F_t(x)$, with
\begin{align}
    x(t) = F_t(x_0) = x_0 + \int_0^t f(x_\tau) d\tau
    ~.
\end{align}

The Koopman operator evolves observable functions $g(x)$ as 

\begin{align}
    \mathcal{K}_t g = g \circ F_t
    ~.
\end{align}

In general, the observables $g$ are elements of an infinite-dimensional function space, typically $L^2$, and so the Koopman operator $\mathcal{K}_t$ is infinite-dimensional. The \emph{extended dynamic mode decomposition} (EDMD) algorithm constructs a finite-dimensional matrix approximation $K$ to the Koopman operator $\mathcal{K}_t$ using a least-squares regression. 

EDMD requires snapshot input/output pairs of data given as a collection of $M$ samples $\{(x_i, y_i)\}$, where $x_i$ is a sampled system state and $y_i$ is its time-shift, $y_i = F_t(x_i)$. The samples are then used to construct the data matrices $\mathbf{X}$ and $\mathbf{Y}$ as 

\begin{align}
    \mathbf{X} = 
    \begin{bmatrix}
        | & | & \; & | \\ 
        x_1 & x_2 & \cdots & x_M \\ 
        | & | & \; & |
    \end{bmatrix}
    ~,
\end{align}

\noindent and 

\begin{align}
    \mathbf{Y} = 
    \begin{bmatrix}
        | & | & \; & | \\ 
        y_1 & y_2 & \cdots & y_M \\ 
        | & | & \; & |
    \end{bmatrix}
    ~,
\end{align}

Next, EDMD requires a \emph{dictionary} $\Psi = [\psi_1 \; \psi_2 \; \cdots \; \psi_N]^\top$ of scalar basis functions $\psi_j: \mathcal{X} \rightarrow \mathbb{R}$. Applying the dictionary to a state vector $x \in \mathcal{X}$ produces the \emph{lifted state} $\Psi(x) = [\psi_1(x) \; \psi_2(x) \; \cdots \; \psi_N(x)]^\top$. Note that $N$ is typically larger than the dimensionality of $x$. Applying the dictionary to the data matrices $\mathbf{X}$ and $\mathbf{Y}$ produces the lifted data matrices 

\begin{align}
    \mathbf{\Psi}_X = 
    \begin{bmatrix}
        | & | & \; & | \\ 
        \Psi(x_1) & \Psi(x_2) & \cdots & \Psi(x_M) \\ 
        | & | & \; & |
    \end{bmatrix}
    ~,
\end{align}

\noindent and 

\begin{align}
    \mathbf{\Psi}_Y = 
    \begin{bmatrix}
        | & | & \; & | \\ 
        \Psi(y_1) & \Psi(y_2) & \cdots & \Psi(y_M) \\ 
        | & | & \; & |
    \end{bmatrix}
    ~,
\end{align}

\noindent with $\mathbf{\Psi}_X$ and $\mathbf{\Psi}_Y$ both $N \times M$ matrices. EDMD finds the best-fit matrix $K$ that evolves the dictionary functions by minimizing the error 

\begin{align}
    &\sum_{i=1}^M || \mathcal{K}\Psi(x_i) - K \Psi(x_i)||^2 \\ 
    = &\sum_{i=1}^M ||\Psi(y_i) - K \Psi(x_i)||^2 \\
    = &|| \mathbf{\Psi}_Y - K \mathbf{\Psi}_X ||^2_F
    ~.
\end{align}

The least-squares solution is given by 

\begin{align}
    K = \mathbf{\Psi}_Y \mathbf{\Psi}_X^\dagger
    ~,
\end{align}

where $\mathbf{\Psi}_X^\dagger$ is the pseudoinverse of $\mathbf{\Psi}_X$. The EDMD matrix $K$ limits to a Galerkin projection of the Koopman operator $\mathcal{K}$ onto the subspace of functions spanned by the dictionary $\Psi$. 

\subsection{EDMD with Control}

Now consider a controlled dynamical system $\dot{x} = f(x, \upsilon)$, where the evolution of the system state $x$ is influenced by an external driver $\upsilon$. This specifies the instantaneous rate of change of the state vector $x$ given an instantaneous control input $\upsilon$. However, the time-evolution of $\upsilon_t$ is essentially arbitrary (as long as it is not unphysical), and it is provided externally. Therefore, flow maps, and by extension Koopman operators, must be defined relative to a given control policy. 
Denote a continuous-time control policy of time duration $T$ as $\vec{\upsilon} := \{\upsilon_t\}_{t \in [0, T]}$. 

The controlled state dynamic $\dot{x} = f(x, \upsilon)$ and initial condition $x_0$, together with a control policy $\vec{\upsilon}$ uniquely defines the state trajectory out to time $T$. Therefore, we can define flow maps parameterized by a control policy $\vec{\upsilon}$ as 

\begin{align}
    F_{\vec{\upsilon}, t}(x) := x + \int_0^t f(x, \upsilon_t) d\tau 
    ~,
\end{align}

\noindent where $\upsilon_t$ is specified from the given control policy $\vec{\upsilon}$. Note this is only valid for times $t \in [0, T]$ within the control horizon. Koopman operators for control systems are thus also parameterized by the given control policy, 

\begin{align}
    \mathcal{K}_{\vec{\upsilon}, t} g := g \circ F_{\vec{\upsilon},t}
    ~,
\end{align}

\noindent where $g(x)$ is an observable function of the state variable only (the parameterized flow map $F_{\vec{\upsilon},t}$ is also a function of just the state variable). In continuous time, the Koopman generator for control systems is parameterized by a single control input. The instantaneous rate of change of observables for that control input is 
\begin{align}
    \dot{g} = L_\upsilon g = \sum_i f_i(x, \upsilon) \frac{\partial g}{\partial x_i}
    ~.
\end{align}

Numerically solving a differential game as a mixed complimentarity problems requires discrete time steps, even though the underlying dynamical system is in continuous time and we would ideally like to model the instantaneous system response to a single control input. From the definition of the Koopman generator in Eq.~\ref{eq:koop_gen}, flow maps over short time intervals $\Delta t$ define Koopman operators that converge to the continuous generator as $\Delta t \rightarrow 0$. Moreover, we will assume that control policies do not vary much within the time interval $\Delta t$ (excepting a possible measure-zero set of discontinuities). 

Therefore, our strategy is to use flow maps over a small interval, $\Delta t = 0.01$ in all of our experiments, with the control input held constant over that interval $\vec{\upsilon} = \{ \upsilon_t = \upsilon^* \; \forall \ t \in [0, \Delta t]\}$, 
\begin{align}
    F_{\upsilon^*, \Delta t}(x) := x + \int_0^{\Delta t} f(x, \upsilon^*) dt
    ~.
\end{align}

Denote the family of associated unit-step Koopman operators simply as, 

\begin{align}
    \mathcal{K}_{\upsilon^*} g := g \circ F_{\upsilon^*, \Delta t}
    ~.
\end{align}

The dependence on $\Delta t$ is given implicitly to avoid notational clutter, $\mathcal{K}_{\upsilon^*} := \mathcal{K}_{\upsilon^*, \Delta t}$, since we are working in the discretized dynamics with fixed $\Delta t$.

The control policies in this choice of discretized dynamics are piecewise constant on the small interval $\Delta t$ by construction. The time step $\Delta t$ is assumed small enough that continuous control policies can be approximated reasonably well. The approach of using a fixed control policy over a short time window is called a zero-order hold. For a given discretized control policy $\{\upsilon^*_1, \upsilon^*_2, \ldots, \upsilon^*_{T-1}\}$, the dynamics of observables is given by the sequential action of the unit-step Koopman operator with the corresponding fixed control $\{\mathcal{K}_{\upsilon^*_1}, \mathcal{K}_{\upsilon^*_2}, \ldots, \mathcal{K}_{\upsilon^*_{T-1}}\}$.



The evolution of system observables $\Psi(x)$ under $F_{\upsilon^*, \Delta t}$ depends on the given control input, as $\upsilon^*$ determines which unit-step Koopman operator $\mathcal{K}_{\upsilon^*}$ acts on $\Psi(x)$. For data-driven approximations, researchers often assume an additive dependence on $\upsilon^*$ to parameterize $\mathcal{K}_{\upsilon^*}$ so that  

\begin{align}
    \Psi(y) = K \Psi(x) + K_\upsilon \upsilon^*
    ~. 
    \label{eq:addsep_EDMDc}
\end{align}

\noindent with $y = F_{\upsilon^*, \Delta t}(x)$. Because the Koopman evolution of observables is linear and this assumes a linear dependence on the controls, we can estimate $K$ and $K_\upsilon$ using the pseudoinverse solution by simply concatenating the control vector $\upsilon^*$ onto the dictionary $\Psi$. 

This formulation of EDMD with control requires data triples $\{(x_i, \upsilon_i, y_i)\}$, with $y = F_{\upsilon, \Delta t}(x, \upsilon)$. Since we have the equations of motion $f(x, \upsilon)$, we create these data triples by uniformly sampling (i.e., equally spaced samples, not uniform random sampling) state vectors $x_i$ and control inputs $\upsilon_i$ over their respective domains and produce the time-shifted state vector $y_i = F_{\upsilon, \Delta t}(x_i, \upsilon_i)$.

Define $\Omega(x, \upsilon) = [\Psi(x) \; \upsilon]^\top$ and the augmented data matrix

\begin{align}
    \mathbf{\Omega} = 
    \begin{bmatrix}
        | & | & \; & | \\ 
        \Omega(x_1, \upsilon_1) & \Omega(x_2, \upsilon_2) & \cdots & \Omega(x_M, \upsilon_M) \\ 
        | & | & \; & |
    \end{bmatrix}
    ~.
\end{align}

Then, we have 

\begin{align}
    \widetilde{K} = \mathbf{\Psi}_Y \mathbf{\Omega}^\dagger
    ~,
\end{align}

\noindent where $\widetilde{K} = [K \; K_\upsilon]$. This is the best-fit (in a least-squares sense) solution to 

\begin{align}
    \Psi(y) \approx \widetilde{K}\Omega(x, \upsilon) = [K \; K_\upsilon] 
    \begin{bmatrix}
        \Psi(x) \\ 
        \upsilon
    \end{bmatrix}
    ~.
\end{align}

The original dynamic mode decomposition with control (DMDc) algorithm uses the identity function for the basis $\Psi$ and computes the pseudoinverse using SVD~\cite{proctor16a}. There, they do two separate SVD rank reductions; one to regularize the pseudoinverse computation, and the other to reduce the dimensionality of the dictionary $\Psi$. The second rank reduction is appropriate, and often necessary, for high-dimensional state vectors $x$. In this work, the state vector is low-dimensional and we explicitly control the size of $\Psi$, so we do not perform this rank reduction. We compute the pseudoinverse by simply calling the \texttt{/} operator in Julia. 

\subsection{Control Dependence for Turret Defense Games}

\begin{figure}
\centering
\includegraphics[width=0.9\textwidth]{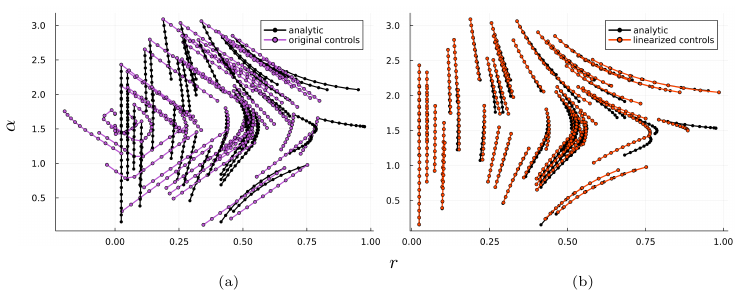}
\caption{Comparison of analytic $r$-$\alpha$ trajectories with EDMD trajectories for (a) the original formulation of the controls, and (b) the linearized formulation with the agent controls split into components.}
\label{fig:EDMD_cntrl_compre}
\end{figure}

For the turret defense game, the state dynamics are largely driven by the controls. In addition, the state dynamics have a highly nonlinear dependence on the controls. Unsurprisingly then, we found that the additive parameterization of $\mathcal{K}_{\upsilon^*}$ in Eq.~\ref{eq:addsep_EDMDc} did not yield an effective approximation of $\mathcal{K}_{\upsilon^*}$ for the turret defense game; see Figure~\ref{fig:EDMD_cntrl_compre}.

To remedy this deficiency, we first removed the nonlinearity of the agent heading $v$ by splitting it into its orthogonal components, 

\begin{align*}
    &v_A \cos v \rightarrow v \\
    &v_A \sin v \rightarrow v_\perp
    ~.
\end{align*}

The control vector is now three-dimensional, $\upsilon = [v \; v_\perp \; u]$. However, EDMD with control using this modified control vector also did not produce an effective approximation of  $\mathcal{K}_{\upsilon^*}$. Although the nonlinearity of the agent heading has been removed, there is a multiplicative relation between the system state and controls in the original equation of motion, whereas the parameterization in Eq.~\ref{eq:addsep_EDMDc} assumes an additive relation between controls and state observables. 

Therefore, we create the following ``lifted'' control vector 
\begin{align}
    \tilde{\upsilon} = 
    \begin{bmatrix}
        \tilde{v} \; \; \;\\
        \tilde{v}_\perp \\
        u \;\;\; 
    \end{bmatrix}
    :=
    \begin{bmatrix}
        r^2 v_A \cos v \\
        r v_A \sin v \\
        u
    \end{bmatrix}
    ~.
\end{align}
The approximation of $\mathcal{K}_{\upsilon^*}$ then has the form
\begin{align}
    \psi(y) = K \psi(x) + K_\upsilon \tilde{\upsilon}
    ~.
\end{align}
Using the analogous augmented data matrix with $\tilde{\Omega}(x, \upsilon) = [\Psi(x) \; \tilde{\upsilon}]^\top$, the best-fit approximation is given by the pseudoinverse solution
\begin{align}
    [K \; K_\upsilon] = \mathbf{\Psi}_Y \tilde{\mathbf{\Omega}}^\dagger
    ~.
\end{align}

Both components $K$ and $K_\upsilon$ now depend on the state vector, but only $K_\upsilon$ depends on the controls. The idea is to isolate the minimal form of dependence on the controls in $K_\upsilon$ required to achieve a reasonable approximation of the unit-step Koopman operator. 
The dependence of the state dynamics on the controls needs to be isolated, and ideally linear. This is the main motivation for a Koopman-based approach in the first place.



\section{Finite-Basis Resolvent Calculations}
\label{app:resolve}

\subsection{Function and Generator Approximations}

With fixed RBF centroids at $x_i$ and fixed evaluations points $x_j$, we have, for any function $g \left(x\right)$,

\begin{align}
    g \left(x\right) &\approx \sum_i g_i \ RBF_i \left(x;x_i\right) \\
    G_{ji} & \equiv RBF_i \left(x_j;x_i\right) \\
    \Rightarrow g\left(x_j\right) & \approx \sum_i G_{ji} g_i
\end{align}

We can also calculate our Koopman generator $L$ and its finite-dimensional approximation $L_{ki}$ this way. If we have dynamics $\dot{x} = f \left(x\right)$ and $L = \sum_k f_k \frac{\partial}{\partial x_k}$, and if we define $\partial G_{jik}$ such that

\begin{gather}
    \frac{\partial g \left(x_j\right)}{\partial x_k} = \sum_i \partial G_{jik} g_i
\end{gather}

\noindent then

\begin{align}
    L g\left(x_j\right) &= \sum_k f_k \left(x_j\right) \frac{\partial g \left(x_j\right)}{\partial x_k} \\
    \Rightarrow L g \left(x_j\right) &= \sum_{i,k} f_k \left(x_j\right) G_{jik} g_i \\
    \Rightarrow \sum_k G_{jk} \left(\sum \limits_i L_{ki} g_i\right) &= \sum_{i,k} f_k \left(x_j\right) \partial G_{jik} g_i
\end{align}

This needs to hold for all functions $g$, so

\begin{align}
    \sum_k G_{jk} L_{ki} &= \sum_k f_k \left(x_j\right) \partial G_{jik}
\end{align}

Assuming that there are a sufficient number of evaluation points $j$, we can solve for $L_{ki}$ as the solution to an overdetermined system of linear equations (multiplying the righthand side by the pseudo-inverse of $G$).

\subsection{Resolvent Integration and Matrix Inversions}

The resolvent calculations include matrix inversions. Our optimization then has two ways that it can handle this. The first option simply performs the calculations ``under the hood'' as part of evaluating the objective function; the necessary derivatives get propagated automatically through the matrix inversions and passed to the optimizer. The matrix inversions are also performed to a level of precision corresponding to Julia's numerical tolerances for such operations in general. This makes the objective function more expensive to evaluate, and it creates a significant amount of nonlinearity, but it avoids adding additional variables or constraints to the problem (thus leaving the size of the problem unchanged). In principle, we could also conduct successive optimizations at progressively greater levels of fidelity, which could potentially speed up the optimization process, without changing the underlying problem.

Alternatively, we could embed the matrix inversions into the optimization as a set of equality constraints. Instead of an objective that looks like

\begin{align}
    \sum_i &c_i A_i^{-1} \left(x\right) b\left(x\right)
\end{align}

\noindent we get a set of constraints and an objective in the respective forms

\begin{align}
    A_i \left(x\right) \gamma_i  &= b_i \left(x\right) \qquad \forall \ i \\
    \sum_i & c_i \gamma_i 
\end{align}

This adds additional constraints and variables ($\gamma_i$) to the problem, and the number added of each will be proportional to the number of integration nodes in the resolvent calculation -- which can increase the size of the problem by orders of magnitude -- but the objective is now linear in those new variables. Furthermore, the tolerance on the equality constraint satisfaction may be worse than that on Julia's linear solve, and changing the fidelity of the resolvent calculation would result in a different problem (because the number of integration nodes would change). That being said, the optimization may be easier to solve, numerically, because the nonlinearity associated with the explicit matrix inversion is no longer present. The optimizer might also be able to violate the constraints, temporarily, to more easily find the optimal solution.

For the purpose of this paper, we focused on the first option. Both methods struggled with convergence and computational cost, but the second option performed more poorly -- especially with the warm start optimizations. The computational cost of the optimization motivated us to use a multi-fidelity approach whereby we used solutions with larger tolerances on both optimality and resolvent error as warm starts for optimizations with tigheter tolerances. This provided strong motivation to use the first option.

\subsection{Alternating Optimizations}

When doing these iterated and alternating optimizations, we copied the primal variables to improve the efficiency of the warm start. These values would be perturbed slightly by IPOPT, and the the barrier function coefficient in IPOPT would also reset at the beginning of a new optimization. The combination of these two things struck a balance between perturbing the warm start solution too much or too little. When alternating back and forth, we would also manipulate the constraints to prevent pathological behaviors. Some constraints were only applicable to one of the players -- in particular, the $\dot{r} \leq 0$ at $r= 1$ constraint could only be affected by the agent's control policy, not the turret's. If the agent's optimization was infeasible with respect to this constraint when it hit its maximum iteration limit, including this constraint in the turret's optimization would make that optimization infeasible regardless of what the turret does. Therefore, the upper and lower bounds on the turret control were essentially rendered inactive during the agent's optimization, and the agent velocity constraint and $r \leq 1$ constraint were essentially rendered inactive during the turret's optimization. For the constraints on $\dot{\alpha}$, during one player's optimization, the other player's control variables would be held constant, and this did not seem to restrict the optimization unduly.

\subsection{Variable-Precision Solution}
\label{sec: variable precision}

In principle, the resolvent-based approach could be used with something like ApproxFun, which would enable more precise error control via dynamic adjustment of the basis size. This could be advantageous in some circumstances, but it would likely not be compatible with a numerical optimizer like IPOPT; EDMD or Deep-DMD, with fixed bases, can make use of these off-the-shelf optimizers. It would instead be necessary to implement a custom optimization algorithm implementation that could work on the relevant functions directly (and thereby avoid the need for a fixed basis size). 

\subsubsection{Derivative Calculations}

A custom optimization algorithm implementation could be done, but it would, in a sense, be re-inventing the wheel. It would then also be necessary to use the resolvent to calculate the ``derivative'' of the cost function with respect to the control function:

\begin{align}
    \dot{x} & = f \left(x,u,v\right) \\
    L &= f_i \left(x,u,v\right) \frac{\partial}{\partial x_i} \\
    \frac{\partial L}{\partial u} &= \frac{\partial f_i}{\partial u} \frac{\partial}{\partial x_i} \\
    \frac{\partial L}{\partial v} &= \frac{\partial f_i}{\partial v} \frac{\partial}{\partial x_i} \\
    F \left(x_0,z\right) & \equiv \left(Iz - L\right)^{-1} \left( g \left(x_0\right) + \frac{1}{z}\left[ 1 - e^{-zT}\right] h \left(x_0\right) \right) \\
    \frac{\partial F \left(x_0,z\right)}{\partial \left(\cdot\right)} &= \left(Iz - L\right)^{-1} \frac{\partial L}{\partial \left(\cdot\right)} F \left(x_0,z\right) \\
    \frac{\partial J \left(x_0,T\right)}{\partial \left(\cdot\right)} &= \frac{1}{2 \pi i} \int_\gamma e^{zT} \frac{\partial F \left(x_0,z\right)}{\partial \left(\cdot\right)} dz \\
    \frac{\partial J \left(T\right)}{\partial \left(\cdot\right)} &= \int_X \left[\frac{1}{2 \pi i} \int_\gamma e^{zT} \frac{\partial F \left(x_0,z\right)}{\partial \left(\cdot\right)} dz\right] dx
\end{align}

Using this to identify ascent and descent directions in the functional space would be relatively straightforward, but extending concepts like Armijo conditions or conjugate directions to an infinite-dimensional functional space would be less obvious. For example, if we want to use a Newton's Method approach, we would need second derivatives:

\begin{align}
    \frac{\partial^2 L}{\partial \left(\cdot\right)} &= \frac{\partial^2 f_i}{\partial \left(\cdot\right)^2} \frac{\partial}{\partial x_i} \\
    \left( g \left(x_0\right) + \frac{1}{z}\left[ 1 - e^{-zT}\right] h \left(x_0\right) \right) & = \left(Iz - L\right) F \left(x_0,z\right)\\
    0 &= - \frac{\partial L}{\partial \left(\cdot\right)} F \left(x_0,z\right) + \left(Iz - L \right) \frac{\partial F \left(x_0,z\right)}{\partial \left(\cdot\right)}\\
    0 &= - \frac{\partial^2 L}{\partial \left(\cdot\right)^2} F \left(x_0,z\right) - 2 \frac{\partial L}{\partial \left(\cdot\right)} \frac{\partial F \left(x_0,z\right)}{\partial \left(\cdot\right)} \nonumber \\
    &+ \left(Iz - L\right) \frac{\partial^2 F \left(x_0,z\right)}{\partial \left(\cdot\right)^2} \\
    \Rightarrow \frac{\partial^2 F \left(x_0,z\right)}{\partial \left(\cdot\right)^2} &= \left(Iz - L\right)^{-1} \left[ \frac{\partial^2 L}{\partial \left(\cdot\right)^2} F \left(x_0,z\right) + 2 \frac{\partial L}{\partial \left(\cdot\right)} \frac{\partial F \left(x_0,z\right)}{\partial \left(\cdot\right)} \right] \\
    \frac{\partial^2 J \left(T\right)}{\partial \left(\cdot\right)^2} &= \int_X \left[\frac{1}{2 \pi i} \int_\gamma e^{zT} \frac{\partial^2 F  \left(x_0,z\right)}{\partial \left(\cdot\right)^2} dz\right] dx
\end{align}

It is not necessarily clear how to use these second derivatives, though. Newton's method in 1-D involves dividing the first derivative by the second derivative. Straight division of functions, however, would create problems if the second derivative of $J$ is zero anywhere. Alternatively, we could consider the gradient as an infinite-dimensional vector and the Hessian as an infinite-dimensional matrix. If we have a control policy $u\left(x\right) = \sum_i u_i \psi_i \left(x\right)$, then we can define the difference between $\frac{\partial }{\partial u}$ and $\frac{\partial }{\partial u_a}$ and define some relevant derivatives as a result:

\begin{align}
    \frac{\partial u}{\partial u_a} &= \psi_a \left(x\right) \\
    \frac{\partial }{\partial u_a} &= \frac{\partial u}{\partial u_a} \frac{\partial }{\partial u}  \nonumber \\
    &= \psi_a \left(x\right) \frac{\partial }{\partial u}  \\
    \frac{\partial L_{ij}}{\partial u} \psi_j \left(x\right) &= \sum_k \frac{\partial \psi_j}{\partial x_k} \frac{\partial f_k \left(x,u\left(x\right),v \left(x\right)\right)}{\partial u} \\
    \frac{\partial f_k \left(x,u\left(x\right)\right)}{\partial u} &= \sum_i \frac{\partial f_{i,k}}{\partial u} \psi_i \left(x\right) \\
    \Rightarrow \frac{\partial L_{ij}}{\partial u} &= \sum_k \frac{\partial \psi_j}{\partial x_k} \frac{\partial f_{i,k}}{\partial u} \\
    \frac{\partial L_{ij}}{\partial u_a} &= \sum_k \frac{\partial \psi_j}{\partial x_k} \frac{\partial f_{i,k}}{\partial u} \psi_a \left(x\right) \nonumber \\
    &= \frac{\partial L_{ij}}{\partial u} \psi_a \left(x\right)
\end{align}

We can then use this to calculate the (infinite-dimensional) gradients and Hessians:

\begin{align}
    \frac{\partial^2 L_{ij}}{\partial u^2} \psi_j \left(x\right) &= \sum_k \frac{\partial \psi_j}{\partial x_k} \frac{\partial^2 f_k \left(x,u\left(x\right),v\left(x\right)\right)}{\partial u^2} \\
    \frac{\partial^2 f_k \left(x,u\left(x\right)\right)}{\partial u^2} &= \sum_i \frac{\partial^2 f_{i,k}}{\partial u^2} \psi_i \left(x\right) \\
    \frac{\partial^2 L_{ij}}{\partial u_a \partial u_b} &= \frac{\partial}{\partial u_b} \left( \frac{\partial L_{ij}}{\partial u_a}\right) \nonumber \\
    &= \frac{\partial u}{\partial u_a} \left( \frac{\partial L_{ij}}{\partial u} \psi_a \left(x\right) \right) \nonumber \\
    &= \frac{\partial L_{ij}}{\partial u^2} \psi_a \left(x\right) \psi_b \left(x\right) \\
    \frac{\partial F \left(x_0,z\right)}{\partial u_a} &= \left(Iz - L\right)^{-1} \frac{\partial L}{\partial u_a} F \left(x_0,z\right) \\
    \frac{\partial^2 F \left(x_0,z\right)}{\partial u_a \partial u_b} &= \left(Iz - L\right)^{-1} \left[ \frac{\partial^2 L}{\partial u_a \partial u_b} F \left(x_0,z\right) + \frac{\partial L}{\partial u_a} \frac{\partial F \left(x_0,z\right)}{\partial u_b} + + \frac{\partial L}{\partial u_b} \frac{\partial F \left(x_0,z\right)}{\partial u_a} \right] \\
    \frac{\partial J \left(T\right)}{\partial u_a} &= \int_X \left[\frac{1}{2 \pi i} \int_\gamma e^{zT} \frac{\partial F \left(x,z\right)}{\partial u_a} dz\right] dx \\
    \frac{\partial^2 J \left(T\right)}{\partial u_a \partial u_b} &= \int_X \left[\frac{1}{2 \pi i} \int_\gamma e^{zT} \frac{\partial^2 F \left(x,z\right)}{\partial u_a \partial u_b} dz \right] dx
\end{align}

Analogous calculations apply to control policy $v$.

\subsubsection{Constraint Implementations}

To enforce constraints on an optimization, there are basically three different approaches: projection, penalty function, and barrier methods. In a finite-dimensional optimization, projection methods are fairly simple: the optimization searches until it hits a constraint boundary, it moves along the constraint boundary by projecting the search direction onto the constraint surface, and it imposes correction steps as need be if the search produces an infeasible point. In the infinite-dimensional case, this is more complicated. Consider a constraint $u \left(x\right) \leq \alpha \left(x\right)$, and let us assume that we have a search direction $\Delta u$. Then, using $R \left(\cdot\right)$ as the ReLU function, we get

\begin{align}
    \Delta c \left(x\right) &= \alpha \left(x\right) - u \left(x\right) \\
    \Delta \hat{u} \left(x\right) &= \Delta u \left(x\right) - R \left( \Delta u \left(x\right) - \Delta c \left(x\right) \right) \nonumber \\
    &= \left\{ \begin{array}{cc} \Delta u \left(x\right) & \Delta u \left(x \right) \leq c \left(x\right) \\ \Delta c \left(x\right) & \Delta u \left(x\right) > c \left(x\right) \end{array} \right.
\end{align}

Note that, in this case, if we are violating the constraint, the $\Delta \hat{u} \left(x\right)$ search direction will push the solution back to a feasible state; we can also implement the ReLU function as 

\begin{gather}
R\left(x\right) = \frac{1}{2} \left(\text{abs} \left(x\right) + x\right)
\end{gather}

Because the ReLU function is not smooth, ApproxFun may have difficulty setting up a basis function representation. A smoothed version of this would use a softplus function:

\begin{align}
    S \left(x;k\right) &= \frac{\log \left( 1 + e^{kx}\right)}{k} \\
    \Delta \hat{u} \left(x\right) &= \Delta u \left(x\right) - S \left(\Delta u \left(x\right) - \Delta c \left(x\right);k\right)
\end{align}

\noindent where $k$ is a shape parameter (governing the sharpness of the approximation). In the more general case,

\begin{align}
    c \left(u \left(x\right),x\right) & \geq 0 \\
    c \left( u \left(x\right) + \Delta u \left(x\right),x\right) & \approx c \left(u \left(x\right),x\right) + \frac{\partial c}{\partial u} \Delta u \\
    \Rightarrow \Delta c \left(x\right) & = -\frac{c \left(u \left(x\right),x\right)}{\frac{\partial c}{\partial u}} \\
    \Delta c \left(x\right) \geq 0 & \Rightarrow \Delta \hat{u} \left(x\right) \leq \Delta c \left(x\right) \\
    \Delta c \left(x\right) \leq 0 & \Rightarrow \Delta \hat{u} \left(x\right) \geq \Delta c \left(x\right)
\end{align}

\noindent and the calculations for $\Delta \hat{u}$ proceed as before. The applications to $v \left(x\right)$ are analogous. If we have known bounds $\Delta c_{min} \left(x\right) \leq \Delta \hat{u} \left(x\right) \leq \Delta c_{max} \left(x\right)$, then we can approximate the hard bounds with a sigmoid function. We also want the sigmoid function to approximate the identity function at the midpoint between $\Delta c_{max}$ and $\Delta c_{min}$. We can do this as follows:

\begin{gather}
    \Delta \hat{u} \left(x\right) = \left[\frac{\Delta c_{max} \left(x\right) - \Delta c_{min} \left(x\right)}{2}\right] \nonumber \\
    \times \text{tanh} \left( \frac{2}{\Delta c_{max} \left(x\right) - \Delta c_{min} \left(x\right)} \left( \Delta u \left(x\right) - \frac{\Delta c_{max} \left(x\right) + \Delta c_{min} \left(x\right)}{2} \right) \right) \nonumber \\ + \frac{\Delta c_{max} \left(x\right) 
    + \Delta c_{min} \left(x\right)}{2}
\end{gather}

Multiple constraints, or multiple $u$ variables, will likely need to be handled on a case-by-case basis.

Alternatively, we could add a (smoothed) penalty function

\begin{align}
    J & \rightarrow J + \rho_{pen} S \left(-c \left(u \left(x\right),x\right)\right) \\
    \frac{\partial J}{\partial u} & \rightarrow \frac{\partial J}{\partial u} - \frac{\rho_{pen}}{1 + e^{-k c \left(u \left(x\right),x\right)}} \frac{\partial c}{\partial u}
\end{align}

\noindent or a barrier function

\begin{align}
    J & \rightarrow J - \rho_{bar} \log \left(c \left(u \left(x\right),x\right)\right) \\
    \frac{\partial J}{\partial u} & \rightarrow \frac{\partial J}{\partial u} - \frac{\rho_{bar}}{c \left(u \left(x\right),x\right)} \frac{\partial c}{\partial u}
\end{align}

These approaches would be simpler and more general, but they would also require some kind of progressive relaxation on the penalty or barrier function weights ($\rho_{pen}$ or $\rho_{bar}$, respectively) to converge to a true optimum of the original problem.

\subsection{Imposing Terminal Conditions}

In the main body of the paper, we considered a fixed-time differential game. This was necessary for both the resolvent- and complementarity-based methods. However, in the future, there may be a way to lift this restriction for the resolvent-based method, and we can illustrate what this might look like via a simple example. Assume that we have a 2D dynamical system with ``cost'' function $g \left(x\right)$

\begin{align}
    \dot{x}_1 & = x_1 \\
    \dot{x}_2 &= x_2 \\
    g \left(x\right) &= x_1 + x_2 \\
    \Rightarrow x_1 \left(t\right) &= x_{1,0} e^t \\
    \Rightarrow x_2 \left(t\right) &= x_{2,0} e^t
\end{align}

\noindent defined on $\left[-1,1\right]^2$ -- so the boundary is $x_1 = \pm 1$, $x_2 = \pm 1$. In this case, the time to the boundary is

\begin{gather}
    t\left(x_0\right) = \max \left(-\ln x_{1,0},-\ln x_{2,0} \right)
\end{gather}

In a basis $\left\{x_1,x_2\right\}$, our Koopman operator and cost function are

\begin{align}
    K_t &= \left[ \begin{array}{cc} e^t &0 \\ 0 &e^t \end{array} \right] \\
    g \left(x\right) &= \left\{ \begin{array}{c} 1 \\ 1\end{array} \right\} \\
    g_t \left(x_0\right) &= \left\{ \begin{array}{c} e^t \\ e^t \end{array} \right\}
\end{align}

What we want to do is make the time of integration dependent on the initial condition -- $t \left(x_0\right)$. If we do this,

\begin{align}
    g_{t \left(x_0\right)} \left(x_0\right) &= \left\{ \begin{array}{c} \max \left(\frac{1}{x_{1,0}},\frac{1}{x_{2,0}}\right) \\ \max \left(\frac{1}{x_{1,0}},\frac{1}{x_{2,0}}\right) \end{array} \right\}
\end{align}

Moving out of the basis, we get

\begin{gather}
    g_{t \left(x_0\right)} \left(x_0\right) = \left\{ \begin{array}{cc} 1 + \frac{x_{1,0}}{x_{2,0}} & x_{2,0} > x_{1,0} \\ 1 + \frac{x_{2,0}}{x_{1,0}} & x_{2,0} < x_{1,0} \end{array} \right\}
\end{gather}

This function is not smooth at $x_{1,0} = x_{2,0}$, but it is continuous there; it will be undefined at the origin -- this approach may require that the domain does not contain a critical point. More generally, this would look like

\begin{align}
    g_{terminal} &= \frac{1}{2 \pi i} \int_\gamma e^{z t_{terminal} \left(x\right)} \left(I z - L \right)^{-1} g \left(x\right) dz \nonumber \\
    &\approx \sum_j \beta_j \left(t_{terminal} \left(x\right)\right) \left(I z_j - L\right)^{-1} g \left(x\right) \\
    g \left(x\right) &= \sum_k g_k \psi_k \left(x\right) \\
    \Rightarrow g_{terminal} \left(x\right) &= \sum_k \left[ \sum_j \beta_j \left(t_{terminal} \left(x\right)\right) \eta_{jk} \right] \psi_k\left(x\right) \\
    \sum_k \left[ \left(I z_j - L \right)^{-1} \right]_{lk} \eta_{jk} & = g_l \
\end{align}

The first challenge would be to come up with an expression for $t_{terminal} \left(x\right)$. The second challenge would be to adapt the numerical integration weights $\alpha_j$ so that they can be variable. This will likely increase the computational costs of the resolvent-based approach significantly, as it adds a lot of additional complexity to the nonlinear optimization. However, there may be contexts in which it is worthwhile to do so.

\bibliography{bibliography}

\end{document}